\documentclass[a4paper,10pt]{article}
\usepackage[a4paper]{geometry}
\usepackage{mathtools}
\usepackage{empheq}
\usepackage{amsfonts,mathrsfs,amssymb,yhmath,stmaryrd}
\usepackage{graphicx}

\newcommand{\ds}{\displaystyle}
\def\nm{\noalign{\medskip}}

\newtheorem{proposition}{Proposition}[section]

\newtheorem{thm}{Theorem}[section]
\newtheorem{rem}{Remark}[section]
\newcommand{\G}[2]{G_{2\omega}^{(0)} (#1,#2)}

\title{Backpropagation Imaging in Nonlinear Harmonic Holography in the Presence of
Measurement and Medium Noises\thanks{\footnotesize This work was
supported  by the ERC Advanced Grant MULTIMOD--267184.}}
\author{Habib Ammari\thanks{\footnotesize Department of Mathematics and Applications,
Ecole Normale Sup\'erieure, 45 Rue d'Ulm, 75005 Paris, France
(habib.ammari@ens.fr, pierre.millien@ens.fr).}
 \and Josselin Garnier\thanks{\footnotesize
Laboratoire de Probabilit\'es et Mod\`eles Al\'eatoires \&
Laboratoire Jacques-Louis Lions, Universit\'e Paris VII, 75205
Paris Cedex 13, France (garnier@math.univ-paris-diderot.fr).} \and
Pierre Millien\footnotemark[2]}

\begin{document}

\maketitle

\begin{abstract}
In this paper, the detection of a small reflector in a randomly
heterogenous medium using second-harmonic generation is
investigated. The medium is illuminated by a time-harmonic plane
wave at frequency $\omega$. It is assumed that the reflector has a
non-zero second-order nonlinear susceptibility, and thus emits a
wave at frequency $2\omega$ in addition to the fundamental
frequency linear scattering. It is shown how the fundamental
frequency signal and the second-harmonic signal propagate in the
medium. A statistical study of the images obtained by migrating
the boundary data is performed. It is proved that the
second-harmonic image is more stable with respect to medium noise
than the one obtained with the fundamental signal. Moreover, the
signal-to-noise ratio for the second-harmonic image does not
depend neither on the second-order susceptibility tensor nor on
the volume of the particle.
\end{abstract}

\bigskip

\noindent {\footnotesize Mathematics Subject Classification
(MSC2000): 35R30, 35B30}


\section{Introduction}

Second-harmonic microscopy is a promising imaging technique based
on a phenomenon called second-harmonic generation (SHG) or
frequency-doubling. SHG requires an intense laser beam passing
through a material with non vanishing second-order susceptibility
\cite{GeneralConsideration}. A second electromagnetic field is
emitted at exactly twice the frequency of the incoming field.
Roughly speaking,
\begin{equation}
\mathbf{E}_{2\omega} \sim \mathbf{E}_\omega \chi^{(2)}
\mathbf{E}_\omega,
\end{equation} where $\chi^{(2)}$ is the second-order susceptibility tensor.
A condition for an object to have non vanishing second-order
susceptibility tensor is to have a noncentrosymmetric structure.
Thus SHG only occurs in a few types of physical bodies: crystals
\cite{miller1964}, interfaces like cell membranes \cite{heinz1991,
campagnola2003}, nanoparticle \cite{zavelani2008, hui2004}, and
natural structures like collagen or neurons
\cite{brown2003,mertz2004}. This makes SHG a very good contrast
mechanism for microscopy, and has been used in biomedical imaging.
SHG signals have a very low intensity because the coefficients in
$\chi^{(2)}$ have a typical size of picometer $/V$
\cite{choy1976}. This is the reason why a high intensity laser
beam is required in order to produce a second-harmonic field that
is large enough to be detected by the microscope. Second-harmonic
microscopy has several advantages. Among others, the fact that the
technique does not involve excitation of molecules so it is not
subject to phototoxicity effect or photobleaching. The excitation
uses near infrared light which has a very good penetration
capacity, and a lot of natural structures (like collagen for
instance) exhibit strong SHG properties, so there is no need for
probes or dyes in certain cases. SHG images can be collected
simultaneously with standard microscopy and
two-photon-excitation-fluorescence microscopy for membrane imaging
(see, for instance, \cite{campagnola2003}).

The coherent nature of the SHG signal allows us to use nonlinear
holography for measuring the complex two-dimensional (amplitude
and phase) SHG signal \cite{hsieh2009three, pu2008harmonic}. The
idea is quite similar to conventional linear holography
\cite{cuche1999digital, schnars1994direct}. A frequency doubling
crystal is used to produce a coherent reference beam at the
second-harmonic frequency, which allows to measure the phase of
the one emitted from the reflector \cite{hsieh2011imaging}.

 On the other hand,
since only the dye/membrane produces the second-harmonic signal,
SHG microscopy allows a precise imaging of the dye/membrane, clear
from any scattering from the surrounding medium, contrary to the
fundamental frequency image, where the signal measured is produced
by both the reflector and the medium. As it will be shown in this
paper, this is the main feature which makes second-harmonic
imaging very efficient when it is not possible to obtain an image
of the medium without the dye in order to filter the medium noise.
In practical situations \cite{hsieh2011imaging}, it is not
possible to get an image without the reflector. The main purpose
of this work is to justify that the second-harmonic generation
acts in such situations as a powerful contrast imaging approach.

More precisely, we study the case of a nanoparticle with non
vanishing second-order susceptibility tensor $\chi^{(2)}$ embedded
in a randomly heterogeneous medium illuminated by an incoming
electromagnetic field at a fixed frequency $\omega$. We give
asymptotic formulas for the electromagnetic field diffracted by
the particle and the medium at the fundamental frequency and at
the second-harmonic frequency. Then we use a backpropagation
algorithm in order to recover the position of the particle from
boundary measurements of the fields. We study the images obtained
by backpropagation both in terms of resolution and stability. In
particular,  we elucidate that the second-harmonic field provides
a more stable image than that from fundamental frequency imaging,
with respect to medium noise, and that the signal-to-noise ratio
for the second-harmonic image does not depend neither on
$\chi^{(2)}$ nor on the volume of the particle. The aforementioned
are the main findings of this study.

The paper is organized as follows. In section \ref{sec2} we
formulate the problem of SHG. In section \ref{sec3}, asymptotic
expansions in terms of the size of the small reflector (the
nanoparticle) of the scattered field at the fundamental frequency
and the second-harmonic generated field are derived. In section
\ref{sec5}, we introduce backpropagation imaging functions for
localizing the point reflector using the scattered field at the
fundamental frequency as well as the second-harmonic field. In
section \ref{sec6}, we perform a stability and resolution analysis
of the backpropagation imaging functions.  We show that the medium
noise affects the stability and resolution of the imaging
functions in different ways. We prove that using the
second-harmonic field renders enhanced stability for the
 reconstructed image.  Our main findings are
delineated by a few numerical examples in section \ref{sec7}. The
paper ends with a short discussion.

\section{Problem formulation} \label{sec2}
Consider a small electric reflector $\Omega_r$ embedded in a
randomly heterogeneous medium in $\mathbb{R}^2$. We assume that
the medium has random fluctuations described by a random process $
\mu$  with Gaussian statistics and mean zero. Furthermore, we
assume that $\mu$ is compactly supported in $\mathbb{R}^2$ and let
$\Omega_\mu := \mbox{supp}(\mu)$. We also assume that the
refractive index of the background homogeneous medium
$\mathbb{R}^2 \setminus \overline{\Omega_\mu}$ is $1$.  The medium
is illuminated by a plane wave at frequency $\omega >0$, intensity
$U_I >0$,  and direction $\theta \in \mathbb{S}^1$:
\begin{equation} U_0 (x)= U_I e^{i \omega \theta
\cdot x},\end{equation} with $\mathbb{S}^1$ being the unit circle.
We assume that the incoming plane wave is polarized in the
transverse magnetic direction.
%The electric field at the fondamental
%frequency is given by
%\begin{equation}
%\mathbb{E}(x,\omega)=\frac{c}{i\omega \varepsilon} (\partial_{x_2}U,-\partial_{x_1} U).
%\end{equation}
The small reflector $\Omega_r$ is in $\Omega_\mu$ and has a
refractive index given by
\begin{equation}
 [\sigma_r-1]\textbf{1}_{\Omega_r}(x),
 \end{equation}
 where $\sigma_r$ is the refractive index contrast of the reflector,
 $\Omega_r$ is compactly supported in $\Omega_\mu$ with volume $|\Omega_r|$,
  and $\textbf{1}_{\Omega_r}$ is the  characteristic function of $\Omega_r$.
 The squared refractive index $n(x)$ in the whole space has then the following form:
\begin{equation}
\frac{1}{n(x)}= \left(1+\mu(x) +
[\sigma_r-1]\textbf{1}_{\Omega_r}(x)\right).
\end{equation}

The scattered field $u_s$ generated by the plane wave satisfies
the Helmholtz equation:

 \begin{equation} \label{sgheq1}
 \left\{
\begin{aligned}
  \nabla\cdot\left(([{\sigma_r}-1]\textbf{1}_{\Omega_r} + \mu + 1) \nabla (u_s + U_0)\right)
  +\omega^2  (u_s + U_0) = 0 \quad \mbox{in } \mathbb{R}^2, \\
 \lim\limits_{\vert x\vert \to \infty} \sqrt{\vert
 x \vert}(\frac{\partial u_s}{\partial \vert x \vert} -i\omega
  u_s) =
 0,
\end{aligned}
\right.
\end{equation}

The point reflector also scatters a second field $v$ at frequency
$2\omega$. The field $v$ satisfies, up to $O(||\mu||^2_{L^\infty
(\Omega_\mu)})$, the following Helmholtz equation
\cite{LightWaves,GeneralConsideration,Soussi}:

 \begin{equation} \label{sgheq2}
 \left\{
\begin{aligned}
  \left(\Delta +   \frac{(2 \omega)^2}{[\sigma_r -1] \textbf{1}_{\Omega_r} + 1} (1- \frac{\mu}{[\sigma_r -1]
   \textbf{1}_{\Omega_r} + 1}
  )\right) v=
  \sum_{k,l= 1,2} \chi_{kl} \partial_{x_k} U \partial_{x_l} U \textbf{1}_{\Omega_r} \quad \mbox{in } \mathbb{R}^2,\\
 \lim\limits_{\vert x\vert \to \infty} \sqrt{\vert x \vert}\left(\frac{\partial v}{\partial \vert x \vert} -
 2 i\omega  v\right) =  0,
\end{aligned}
\right.
\end{equation}
where $\chi$ is the electric polarization of the reflector, and
can be written as $\chi(x)=(\chi_{ij})_{i,j=
1,2}\mathbf{1}_{\mbox{r}}(x)$ and $U=u_s+U_0$ is the total field.
Here the second-harmonic field is assumed to be in the transverse
electric mode. The polarization of the second-harmonic field is
given by symmetry properties of the second-order susceptibility
tensor $\chi$. This transverse magnetic--transverse electric
polarization mode is known to be supported by a large class of
optical nonlinear materials
 \cite{shen1984principles}. We choose this polarization mode so that a two-dimensional study of the second
  harmonic generation with scalar fields would be possible. The results would be pretty similar in a general
   three-dimensional case, but the computations would be much elusive.
   The coupled problems (\ref{sgheq1}) and (\ref{sgheq2})
have been mathematically investigated in \cite{bao1, bao3, bao2}.

Let us consider $\Omega$ to be a domain large enough so that
$\Omega_\mu = \mbox{supp}(\mu) \Subset \Omega$ and  measure the
fields $u_s$ and $v$ on its boundary $\partial \Omega$. The goal
of the imaging problem is to locate the reflector from the
far-field measurements of the scattered field $u_s$ at the
fundamental frequency  and/or the second-harmonic generated field
$v$. It will be shown in this paper that the use of the
second-harmonic field yields a better stability properties than
the use of the scattered field at the fundamental frequency in the
presence of medium noise.

\section{Small-volume expansions} \label{sec3}

In this section, we establish small-volume expansions for the
solutions of  problems (\ref{sgheq1}) and (\ref{sgheq2}). We
assume that the reflector is of the form $\Omega_r =z_r+\delta B$,
where its characteristic size $\delta$ is small, $z_r$ is its
location, and $B$ is a smooth domain such that $B\subset B(0,1)$.

\subsection{Fundamental frequency problem}
Let  $U^{(\mu)} =u_s^{(\mu)}+U_0$ be the total field that would be
observed in the absence of any reflector. The scattered field
$u_s^{(\mu)}$ satisfies
\begin{equation} \label{eqdepart} \left\{
\begin{aligned}\nabla\cdot \left((1+\mu) \nabla (u^{(\mu)}_s
+U_0) \right) +  \omega^2 (u^{(\mu)}_s  + U_0) = 0  \quad \mbox{in
}
\mathbb{R}^2, \\
\lim\limits_{\vert x\vert \to \infty} \sqrt{\vert
 x \vert}(\frac{\partial u^{(\mu)}_s}{\partial \vert x \vert} -i\omega
  u^{(\mu)}_s) =
 0.
\end{aligned}
\right.
\end{equation}
Therefore,
$$
\nabla\cdot (1+\mu) \nabla u^{(\mu)}_s  + \omega^2 u^{(\mu)}_s = -
\nabla \cdot \mu \nabla U_0 \quad \mbox{in } \mathbb{R}^2.
$$
Since $\Omega_\mu \Subset \Omega$, the following estimate holds
\begin{equation} \label{estimateborn}
|| u_s^{(\mu)} ||_{H^1(\Omega)}
 \leq C ||\mu||_{L^\infty}
\end{equation}
for some positive constant $C$ independent of $\mu$. Here,
$H^1(\Omega)$ is the set of functions in $L^2(\Omega)$, whose weak
derivatives are in $L^2(\Omega)$. We refer the reader to Appendix
\ref{appenda} for a proof of (\ref{estimateborn}), which uses the
same arguments as those in \cite{abboud2,abboud1}. Actually, one
can prove that
$$
u_s^{(\mu)}(x) = - \int_{\Omega_\mu} \mu(y) \nabla U_0(y) \cdot
\nabla G^{(0)}_\omega(x,y) dy + O (||\mu||_{L^\infty}^2), \quad x
\in \Omega.
$$
Moreover, writing
$$
\nabla\cdot \left((1+\mu) \nabla (u^{(\mu)}_s +U_0) \right) = -
\omega^2 (u^{(\mu)}_s  + U_0),
$$
it follows by using Meyers' theorem \cite{meyers} (see also
\cite[pp. 35-45]{meyers2}) that there exists $\eta>0$ such that
for all $0\leq \eta^\prime \leq \eta$,
$$
\begin{array}{lll}
|| \nabla u_s^{(\mu)} ||_{L^{2+\eta^\prime}(\Omega^\prime)}
 &\leq &|| \nabla (u_s^{(\mu)} + U_0) ||_{L^{2+\eta^\prime}(\Omega)}
 + || \nabla U_0 ||_{L^{2+\eta^\prime}(\Omega)}\\ \nm &\leq& C
 || u_s^{(\mu)} + U_0 ||_{L^{2+\eta^\prime}(\Omega)}
 + || \nabla U_0 ||_{L^{2+\eta^\prime}(\Omega)} \\ \nm &\leq& C
 || u_s^{(\mu)}||_{L^{2+\eta^\prime}(\Omega)} + C^\prime
\end{array}
$$
for some positive constants $C$ and $C^\prime$, where
$\Omega^\prime \Subset \Omega$. From the continuous embedding of
$H^1(\Omega)$ into $L^{2+\eta^\prime}(\Omega)$ and
(\ref{estimateborn}) we obtain
$$
|| u_s^{(\mu)}||_{L^{2+\eta^\prime}(\Omega)} \leq
C^{\prime\prime},
$$
for some constant $C^{\prime\prime}$ independent of $\mu$.
Therefore,
\begin{equation} \label{estimateborn2}
|| \nabla u_s^{(\mu)} ||_{L^{2+\eta^\prime}(\Omega^\prime)}
 \leq C
\end{equation}
for some constant $C$ independent of $\mu$.

Now, on one hand, by subtracting (\ref{sgheq1}) from
(\ref{eqdepart}), we get
\begin{multline} \label{subtract}
\nabla\cdot\left(([{\sigma_r}-1] \textbf{1}_{\Omega_r} + \mu + 1)
\nabla (u_s - u_s^{(\mu)})\right)
  +\omega^2  (u_s - u_s^{(\mu)}) = - \nabla \cdot
  [{\sigma_r}-1] \textbf{1}_{\Omega_r} \nabla U_0  \\ - \nabla \cdot
  [{\sigma_r}-1] \textbf{1}_{\Omega_r} \nabla u_s^{(\mu)}  \quad \mbox{in } \mathbb{R}^2.
  \end{multline}
On the other hand, we have
$$\begin{array}{lll}
|| [{\sigma_r}-1] \textbf{1}_{\Omega_r} \nabla u_s^{(\mu)}
||_{L^2(\Omega)} &\leq &\ds C | \Omega_r|^{\frac{\eta}{8 + 2\eta}}
|| \nabla u_s^{(\mu)} ||_{L^{2+
\frac{\eta}{2}}(\Omega)} \\
\nm &\leq & C | \Omega_r|^{\frac{\eta}{8 + 2\eta}} || \nabla
u_s^{(\mu)} ||_{L^{2}(\Omega)}^{\frac{1}{4+\eta}} || \nabla
u_s^{(\mu)} ||_{L^{2+\eta}(\Omega)}^{\frac{1}{4+\eta}},\\
\end{array}
$$
and hence, by (\ref{estimateborn}) and (\ref{estimateborn2}), we
arrive at
$$
|| [{\sigma_r}-1]\textbf{1}_{\Omega_r} \nabla u_s^{(\mu)}
||_{L^2(\Omega)} \leq C | \Omega_r|^{\frac{\eta}{8 + 2 \eta}}
||\mu||_{L^\infty}^{\frac{2}{4+\eta}}.$$ Therefore, we can neglect
in (\ref{subtract}) the term $\nabla \cdot
[{\sigma_r}-1]\textbf{1}_{\Omega_r} \nabla u_s^{(\mu)}$ as
$||\mu||_{L^\infty} \rightarrow 0$.

Let $w^{(\mu)}$ be defined by
\begin{equation}\label{eqwmu}
\begin{array}{l} \ds \nabla \cdot (1 + \mu+ [\sigma_r -1] \textbf{1}_{\Omega_r}) \nabla w^{(\mu)}  + \omega^2 w^{(\mu)} =
\nabla \cdot [{\sigma_r}-1]\textbf{1}_{\Omega_r} \nabla (x-z_r)
 \quad \mbox{in }
\mathbb{R}^2,
\end{array}
\end{equation}
subject to the Sommerfeld radiation condition. Using the Taylor
expansion
$$
U_0(x) = U_0(z_r) + (x-z_r)\cdot \nabla U_0(z_r) + O(|x-z_r|^2),
$$
one can derive the inner expansion
\begin{equation}
\label{inner} (u_s-u_s^{(\mu)})(x) =  w^{(\mu)}(x) \cdot \nabla
U_0(z_r) + O(\delta^2),
\end{equation}
for $x$ near $z_r$. The following estimate holds. We refer the
reader to Appendix \ref{appendb} for its proof.
\begin{proposition} \label{propb} There exists a positive constant $C$
independent of $\delta$ such that
$$|| u_s - u_s^{(\mu)} -  w^{(\mu)}(x)
\cdot \nabla U_0(z_r) ||_{H^1(\Omega)} \leq C \delta^2.
$$
\end{proposition}

Let  $G_{\omega}^{(\mu)}$ be the outgoing Green function in the
random medium, that is, the solution to
\begin{equation} \label{gmuomega} (\nabla \cdot (1+ \mu) \nabla +\omega^2) G_{\omega}^{(\mu)}(.,z)=-\delta_z \
\quad \mbox{in }  \mathbb{R}^2,
\end{equation}
subject to the Sommerfeld radiation condition. Here, $\delta_z$ is
the Dirac mass at $z$. An important property satisfied by
$G_{\omega}^{(\mu)}$ is the  reciprocity property
\cite{ammarimethods}:
\begin{equation}
\label{reciprocity} G_{\omega}^{(\mu)}(x,z)  =
G_{\omega}^{(\mu)}(z,x), \qquad x\neq z.
\end{equation}

Let us denote by $G_{\omega}^{(0)}$ the outgoing background Green
function, that is, the solution to
\begin{equation} \label{gomega} (\Delta+\omega^2) G_{\omega}^{(0)}(.,z)=-\delta_z \
\qquad \mbox{in } \mathbb{R}^2,
\end{equation}
subject to the Sommerfeld radiation condition.

The Lippmann-Schwinger representation formula:
$$\begin{array}{lll}
(G_\omega^{(\mu)} - G_{\omega}^{(0)})(x,z_r) &=&\ds
\int_{\Omega_\mu} \mu(y) \nabla G_\omega^{(\mu)}(y,z_r)\cdot
\nabla G_{\omega}^{(0)}(x,y)\, dy\\
\nm &=&\ds \int_{\Omega_\mu} \mu(y) \nabla
G_\omega^{(0)}(y,z_r)\cdot \nabla G_{\omega}^{(0)}(x,y)\, dy \\
\nm && \ds + \int_{\Omega_\mu} \mu(y) \nabla (G_\omega^{(\mu)} -
G_{\omega}^{(0)})(y,z_r)\cdot \nabla G_{\omega}^{(0)}(x,y)\, dy
\end{array}
$$
holds for $x \in \partial \Omega$. Since $\Omega_\mu \Subset
\Omega$, we have
$$\begin{array}{l}
\ds \bigg|(G_\omega^{(\mu)} - G_{\omega}^{(0)})(x,z_r) -
\int_{\Omega_\mu} \mu(y) \nabla G_\omega^{(0)}(y,z_r)\cdot \nabla
G_{\omega}^{(0)}(x,y)\, dy \bigg| \leq\\ \nm \qquad\ds
||\mu||_{L^\infty} || \nabla
G_{\omega}^{(0)}(x,\cdot)||_{L^\infty(\Omega_\mu)} || \nabla
(G_\omega^{(\mu)} -
G_{\omega}^{(0)})(\cdot,z_r)||_{L^2(\Omega_\mu)}.
\end{array}
$$
Similarly to (\ref{estimateborn}), one can prove that
\begin{equation} \label{estimateborng} ||
\nabla (G_\omega^{(\mu)} - G_{\omega}^{(0)})(\cdot,z_r)
||_{L^2(\Omega_\mu)}
 \leq C ||\mu||_{L^\infty},
\end{equation}
and hence, there exists a positive constant $C$ independent of
$\mu$ such that
\begin{equation} \label{guome}
\ds \bigg|(G_\omega^{(\mu)} - G_{\omega}^{(0)})(x,z_r) -
\int_{\Omega_\mu} \mu(y) \nabla G_\omega^{(0)}(y,z_r)\cdot \nabla
G_{\omega}^{(0)}(x,y)\, dy \bigg| \leq C ||\mu||_{L^\infty}^2,
\end{equation}
uniformly in $x \in \partial \Omega$.

Since \begin{equation} \label{nabla2} || \nabla \nabla
G_{\omega}^{(0)}(x,\cdot)||_{L^\infty(\Omega_\mu)} \leq C
\end{equation} uniformly in $x \in \partial \Omega$,
the estimate
\begin{equation} \label{guomenabla}
\ds \bigg| \nabla (G_\omega^{(\mu)} - G_{\omega}^{(0)})(x,z_r) -
\nabla \int_{\Omega_\mu} \mu(y) \nabla G_\omega^{(0)}(y,z_r)\cdot
\nabla G_{\omega}^{(0)}(x,y)\, dy \bigg| \leq C
||\mu||_{L^\infty}^2,
\end{equation}
holds in exactly the same way as in (\ref{guome}). Therefore, the
following Born approximation holds.
\begin{proposition} We have
$$\begin{array}{lll}
\ds G_\omega^{(\mu)}(x,z_r) &=&  \ds G_{\omega}^{(0)}(x,z_r) -
\int_{\Omega_\mu} \mu(y) \nabla G_\omega^{(0)}(y,z_r)\cdot \nabla
G_{\omega}^{(0)}(x,y)\, dy  + O( ||\mu||_{L^\infty}^2),\\
\nm \ds \nabla G_\omega^{(\mu)} (x,z_r) &=&\ds \nabla
G_{\omega}^{(0)}(x,z_r) - \nabla \int_{\Omega_\mu} \mu(y) \nabla
G_\omega^{(0)}(y,z_r)\cdot \nabla G_{\omega}^{(0)}(x,y)\, dy +
O(||\mu||_{L^\infty}^2)
\end{array}
$$
uniformly in $x \in \partial \Omega$.
\end{proposition}

We now turn to an approximation formula for $w^{(\mu)}$ as
$||\mu||_{L^\infty} \rightarrow 0$. By integrating by parts we get
$$
w^{(\mu)}(x) = (1-\sigma_r) \int_{\Omega_r} \nabla (w^{(\mu)}(y) -
(y-z_r))\cdot \nabla G_{\omega}^{(\mu)}(x,y)\, dy, \quad x \in
\mathbb{R}^2.
$$
Using (\ref{nabla2}) we have, for $x$ away from $\Omega_r$,
\begin{equation} \label{wmuapp3}
w^{(\mu)}(x) = (1-\sigma_r) [\int_{\Omega_r} \nabla (w^{(\mu)}(y)
- (y-z_r))\, dy] \cdot [\nabla G_{\omega}^{(\mu)}(x,z_r) +
O(\delta)].
\end{equation}

Now let $\textbf{1}_B$ denote the  characteristic function of $B$.
Let $\widetilde{w}$ be the solution to
\begin{equation} \label{fv1conduc}\left\{
\begin{array}{l} \ds \nabla \cdot (1 + [\sigma_r -1] \textbf{1}_B) \nabla  \widetilde{w} = 0 \quad \mbox{in }
\mathbb{R}^2,\\
\nm
 \widetilde{w}(\widetilde{x}) - \widetilde{x}   \rightarrow 0 \quad \mbox{as } |\xi| \rightarrow
+
\infty.\\
\end{array} \right. \end{equation}

The following result holds. We refer the reader to Appendix
\ref{appendc} for its proof.
\begin{proposition} \label{propappendc} We have
\begin{equation}
\label{wmuapp}\nabla\left(  w^{(\mu)}(y) - (y-z_r) \right)= \delta
\nabla \widetilde{w}(\widetilde{y}) + O( \delta [||\mu||_{L^\infty} +
(\delta \omega)^2]),
\end{equation}
where the scaled variable $$\widetilde{y} =
\frac{y-z_r}{\delta}.$$
\end{proposition}

From (\ref{wmuapp}), it follows that
\begin{equation}
\label{wmuapp2} \int_{\Omega_r} \nabla (w^{(\mu)}(y) - (y-z_r))\,
dy = \delta^2 \int_B \nabla \widetilde{w}(\widetilde{x}) \,
d\widetilde{x} + O(\delta^3[ ||\mu||_{L^\infty} + (\delta
\omega)^2]).
\end{equation}

Define the polarization tensor associated to $\sigma_r$ and $B$ by
(see \cite{ammari2004reconstruction})
$$
M(\sigma_r, B) := (\sigma_r -1) \int_B \nabla
\widetilde{w}(\widetilde{x})\, d\widetilde{x} ,
$$
where $\widetilde{w}$ is the solution to (\ref{fv1conduc}). The
matrix $ M(\sigma_r, B)$ is symmetric definite (positive if
$\sigma_r>1$ and negative if $\sigma_r<1$). Moreover, if $B$ is a
disk,  then $M(\sigma_r, B)$ takes the form
\cite{ammari2004reconstruction}:
$$M(\sigma_r, B) =
 \frac{2 (\sigma_r -1)}{\sigma_r+1} \vert B\vert I_2,$$ where $I_2$ is the
 identity matrix.

To obtain an asymptotic expansion of $u_s(x) - u_s^{(\mu)}(x)$ in
terms of the characteristic size $\delta$ of the scatterer, we
take the far-field expansion of (\ref{inner}). Plugging formula
(\ref{wmuapp2}) into (\ref{wmuapp3}), we obtain the following
small-volume asymptotic expansion.
\begin{proposition} \label{propappendb}
%Assume that $\delta \omega \ll 1$ and $||\mu||_{L^\infty(\Omega_\mu)} \ll 1$.
We have
\begin{equation}\label{DAusnf}
u_s(x)=u_s^{(\mu)}(x) - \delta^2 M(\sigma_r, B) \nabla
U_0(z_r)\cdot \nabla G_\omega^{(\mu)} (x,z_r) +O(\delta^3[1 +
||\mu||_{L^\infty} + (\delta \omega)^2]),
\end{equation}
uniformly in $x \in \partial \Omega$.
\end{proposition}

Finally, using (\ref{guomenabla}) we arrive at the following
result.
\begin{thm} \label{propappendbf}
%For $||\mu||_{L^\infty(\Omega_\mu)}$ small enough,
We have as $\delta$ goes to zero
\begin{equation}\label{DAus} \begin{array}{lll}
(u_s-u_s^{(\mu)})(x) &=& - \ds \delta^2 M(\sigma_r, B) \nabla
U_0(z_r)\cdot \bigg[ \nabla G_\omega^{(0)} (x,z_r) + \nabla
\int_{\Omega_\mu} \mu(y) \nabla G_\omega^{(0)}(y,z_r)\cdot \nabla
G_{\omega}^{(0)}(x,y)\, dy  \bigg] \\ \nm && +O(\delta^3 [1 +
||\mu||_{L^\infty} + (\delta \omega)^2] + \delta^2
||\mu||^2_{L^\infty} ),
\end{array}
\end{equation}
uniformly in $x \in \partial \Omega$.
\end{thm}
Theorem \ref{propappendbf} shows that the asymptotic expansion
(\ref{DAus}) is uniform with respect to $\omega$ and $\mu$,
provided that $\omega \leq C/\delta$ and $||\mu||_{L^\infty} \leq
C^\prime \sqrt{\delta}$ for two positive constants $C$ and
$C^\prime$.

\subsection{Second-harmonic problem}

We apply similar arguments to derive a small-volume expansion for
the second-harmonic field at frequency $2 \omega$.

Introduce $G_{2\omega}^{(\sigma_r,\mu)}(.,z)$ the outgoing
solution of $$ \left(\Delta +   \frac{(2 \omega)^2}{[\sigma_r -1]
\textbf{1}_{\Omega_r} + 1} (1- \frac{\mu}{[\sigma_r -1]
\textbf{1}_{\Omega_r} + 1} )\right)
G_{2\omega}^{(\sigma_r,\mu)}(.,z) = -\delta_z \qquad \mbox{in }
\mathbb{R}^2.$$ Let $G_{2\omega}^{(0)}$ be the outgoing solution
to (\ref{gomega}) with $\omega$ replaced by $2\omega$.

Similarly to (\ref{DAus}), an asymptotic expansion for
$G_{2\omega}^{(\sigma_r,\mu)}$ in terms of $\delta$ can be
derived. We have
$$
(G_{2\omega}^{(\sigma_r,\mu)} - G_{2\omega}^{(\mu)}) (x,z) =
O(\delta^2)$$ for $x\neq z$ and $x,z$ away from $z_r$. Here
$G_{2\omega}^{(\mu)}$ is the solution to (\ref{gmuomega}) with
$\omega$ replaced by $2\omega$. Moreover, the Born approximation
yields
$$\begin{array}{lll} (G_{2\omega}^{(\sigma_r,\mu)} -
G_{2\omega}^{(0)})(x,z)&=&\ds - (2\omega)^2 \int_{\Omega_\mu}
\mu(y)G_{2\omega}^{(0)}(y,z)G_{2\omega}^{(0)}(x,y) dy + O(\delta^2
+ ||\mu||_{L^\infty}^2)
\end{array}
$$
for $x\neq z$ and $x,z$ away from $z_r$. From the integral
representation formula:
$$
v(x)= -\int_{\Omega_r}
  \sum_{k,l= 1,2} \chi_{kl} \partial_{x_k} U(y) \partial_{x_l} U(y) G_{2\omega}^{(\sigma_r,\mu)}(x,y) dy,
$$
it follows that
\begin{equation}
v(x)=-\delta^2 |B|  \left(\sum_{k,l} \chi_{kl}
\partial_{x_k} U(z_r) \partial_{x_l} U(z_r) \right) G^{(\sigma_r,\mu)}_{2\omega}(x,z_r) +
O(\delta^3),
\end{equation}
where $|B|$ denotes the volume of $B$, and hence, keeping only the
terms of first-order in $\mu$ and of second-order in
$\delta$:\begin{multline}
 v(x)=-\delta^2 |B| \left(\sum_{k,l}
 \chi_{kl} \partial_{x_k}U(z_r) \partial_{x_l} U(z_r)\right)
  \\ \left[ \G{x}{z_r} - 4\omega^2 \int_\Omega \mu(y) G_{2\omega}^{(0)} (x,y)
  \G{y}{z_r}dy
  + O(||\mu||^2_{L^\infty}) \right]+
  O(\delta^3).
\end{multline}
We denote by $(S)^\theta$ the source term (the source term
strongly depends on the angle $\theta$ of the incoming plane
wave):
\begin{equation}(S)^\theta=\left(\sum_{k,l}\chi_{kl} \partial_{x_k}U(z_r)
\partial_{x_l} U(z_r)\right).
\end{equation}
Now, since
\begin{equation}
U(x)=U_Ie^{i\omega\theta\cdot x} + \int_\Omega \mu(y) \nabla
G_\omega^{(0)}(x,y) \cdot \nabla U_0(y) dy +
O(||\mu||^2_{L^\infty} + \delta),
\end{equation}
which follows by using the Born approximation and the inner
expansion (\ref{inner}), we can give an expression for the partial
derivatives of $U$. We have
\begin{equation}
\partial_{x_k}U(x)=i\omega \theta_k U_I e^{i\omega \theta
 \cdot x} - i \omega \theta \cdot \int_\Omega \nabla (\mu(y) e^{i
 \omega \theta \cdot y}) \partial_{x_k}
G_\omega^{(0)}(x,y)  dy + O(||\mu||^2_{L^\infty} + \delta).
\end{equation}
We can rewrite the source term as
\begin{multline}
\left(\sum_{k,l}\chi_{k,l} \partial_{x_k}U(z_r)\partial_{x_l} U(z_r)
 \right)= - \omega^2 U_I^2 \sum_{k,l} \chi_{kl} \bigg[  \theta_k \theta_l
 e^{i\omega \theta\cdot z_r} \\
- \theta_k \theta \cdot \int_\Omega \nabla (\mu(y) e^{i
 \omega \theta \cdot y}) \partial_{x_l}
G_\omega^{(0)}(z_r,y)  dy - \theta_l \theta \cdot \int_\Omega
\nabla (\mu(y) e^{i
 \omega \theta \cdot y}) \partial_{x_k}
G_\omega^{(0)}(z_r,y)  dy \\
+ \theta \cdot \int_\Omega \nabla (\mu(y) e^{i
 \omega \theta \cdot y}) \partial_{x_l}
G_\omega^{(0)}(z_r,y)  dy \theta \cdot \int_\Omega \nabla (\mu(y)
e^{i
 \omega \theta \cdot y}) \partial_{x_k}
G_\omega^{(0)}(z_r,y)  dy \bigg]
\\
+ O(||\mu||^2_{L^\infty} + \delta).
\end{multline}
Assume that $\mu \in \mathcal{C}^{0,\alpha}$ for $0<\alpha<{1/2}$.
From
\begin{equation}\label{transfoSrand}\begin{array}{lll}
\ds \int_\Omega \nabla (\mu(y) e^{i
 \omega \theta \cdot y}) \partial_{x_l}
G_\omega^{(0)}(z_r,y)  dy  &=&\ds \int_\Omega \nabla (\mu(y) e^{i
 \omega \theta \cdot y} - \mu(z_r) e^{i
 \omega \theta \cdot z_r}) \partial_{x_l}
G_\omega^{(0)}(z_r,y)  dy \\ \nm &=&\ds - \int_\Omega \nabla
\partial_{x_l} G_\omega^{(0)}(z_r,y) (\mu(y) e^{i
 \omega \theta \cdot y} - \mu(z_r) e^{i
 \omega \theta \cdot z_r}) dy
 \end{array}
\end{equation}
one can show that, for $0<\alpha^\prime \leq \alpha$, we have
\cite{trudinger}
$$
\bigg| \theta \cdot \int_\Omega \nabla (\mu(y) e^{i
 \omega \theta \cdot y}) \partial_{x_l}
G_\omega^{(0)}(z_r,y)  dy \theta \cdot \int_\Omega \nabla (\mu(y)
e^{i
 \omega \theta \cdot y}) \partial_{x_k}
G_\omega^{(0)}(z_r,y)  dy  \bigg| \leq C
||\mu||^2_{\mathcal{C}^{0,\alpha^\prime}},$$ where  $C$ is a
positive constant independent of $\mu$.

So, if we split $(S)^\theta$ into a deterministic part and a
random part:
$$(S)^\theta=(S)_{det}^\theta + (S)_{rand}^\theta + O(||\mu||^2_{\mathcal{C}^{0,\alpha}}
 + \delta),$$ we get
\begin{equation}\label{DEFSDET}
 (S)_{det}^\theta = -  \omega^2 U_I^2  e^{i2\omega \theta \cdot z_r} \sum_{k,l} \chi_{k,l} \theta_k
 \theta_l,
\end{equation}
and
\begin{equation} \label{DEFSRAND}
\begin{array}{lll} \ds (S)_{rand}^\theta &= & \ds\omega^2 \sum_{k,l} \chi_{k,l} \bigg[
\theta_k \theta \cdot \int_\Omega \nabla (\mu(y) e^{i
 \omega \theta \cdot y}) \partial_{x_l}
G_\omega^{(0)}(z_r,y)  dy \\ \nm && \ds + \theta_l \theta \cdot
\int_\Omega \nabla (\mu(y) e^{i
 \omega \theta \cdot y}) \partial_{x_k}
G_\omega^{(0)}(z_r,y)  dy\bigg].
\end{array}
\end{equation}
 Finally, we obtain the following result.
\begin{thm} Assume that $\mu \in \mathcal{C}^{0,\alpha}$ for $0<\alpha<{1/2}$.
Let $0<\alpha^\prime\leq \alpha$.  The following asymptotic
expansion holds for $v$ as $\delta$ goes to zero:
\begin{multline}\label{DAv}
 v(x) = -\delta^2 |B| \bigg( (S)_{det}^\theta \left[ \G{x}{z_r}
 - 4\omega^2 \int_\Omega  \mu(y) G_{2\omega}^{(0)} (x,y)
 \G{y}{z_r}dy  \right] + (S)_{rand}^\theta \G{x}{z_r}
 \bigg)\\
 + O(\delta ^3 + \delta^2 ||\mu||^2_{\mathcal{C}^{0,\alpha^\prime}})
\end{multline}
 uniformly in $x\in \partial \Omega$.
\end{thm}

\section{Imaging functional} \label{sec5}

In this section, two imaging functionals are presented for
locating small reflectors. For the sake of simplicity, we assume
that $B$ and $\Omega$ are  disks centered at $0$ with radius $1$
and $R$, respectively.

\subsection{The fundamental frequency case}
We assume that we are in possession of the following data: $\{
u_s(x), \ x\in \partial\Omega\}$. We introduce the
reverse-time imaging functional\begin{equation}\label{DefI}
\forall z^S \in \Omega,\ I(z^S)=\int_{\partial \Omega \times
\mathbb{S}^1} \frac{1}{i\omega} e^{-i\omega \theta \cdot
z^S}\theta^\top \overline{\nabla G_\omega^{(0)} (x,z^S)}
u_s(x) d\sigma(x) d\sigma(\theta),
\end{equation}
where $\top$ denotes the  transpose. Introduce the matrix:
\begin{equation}\label{DEFR}
R_\omega(z_1,z_2)= \int_{\partial \Omega }\overline{\nabla
G_\omega^{(0)} (x,z_1) }  \nabla G_\omega^{(0)} (x,z_2)^\top
d\sigma(x), \qquad z_1,z_2 \in \Omega^\prime \Subset \Omega.
\end{equation}
Using (\ref{DAus}), we have the following expansion for $I(z^S),
z^S \in \Omega^\prime$,

\begin{multline}\label{DAI}
I(z^S)= \int_{\partial \Omega \times
\mathbb{S}^1} \frac{1}{i\omega} e^{-i\omega \theta \cdot
z^S}\theta^\top \overline{\nabla G_\omega^{(0)} (x,z^S)}
u_s^{(\mu)}(x) d\sigma(x) d\sigma(\theta) \\- \frac{2 \pi \delta^2 (\sigma_r-1)}{\sigma_r+1} U_I
\int_{\mathbb{S}^1} e^{-i\omega \theta \cdot (z^S-z_r) }
\theta^\top \bigg[ R_\omega(z^S,z_r)\\
\nm  + \int_{\partial \Omega} \overline{\nabla G_\omega^{(0)}
(x,z^S)} \left( \nabla \int_{\Omega_\mu} \mu (y) \nabla
G_\omega^{(0)}(y,z_r) \cdot  \nabla G_\omega^{(0)}(x,y) dy\right)^\top
d\sigma(x)
\bigg] \theta d\sigma(\theta) \\
+ O(\delta^3 + \delta^2 ||\mu||_{L^\infty}^2).
\end{multline}

Note that
$$\begin{array}{l}\ds
\int_{\partial \Omega} \overline{\nabla G_\omega^{(0)} (x,z^S)}
\left(\nabla \int_{\Omega_\mu} \mu (y) \nabla
G_\omega^{(0)}(y,z_r) \cdot  \nabla G_\omega^{(0)}(x,y) dy \right)^\top
d\sigma(x) \\ \qquad \ds  =   \int_{\Omega_\mu} \mu (y)
\int_{\partial \Omega} \overline{\nabla G_\omega^{(0)}
(x,z^S)} \left(\nabla\nabla G_\omega^{(0)}(x,y) \nabla
G_\omega^{(0)}(y,z_r)\right)^\top d\sigma(x) dy. \end{array}$$

\begin{rem}
Here, the fact that not only we backpropagate the boundary data
but also we average it over all the possible illumination angles
in $\mathbb{S}^1$ has two motivations. As will be shown later in
section \ref{sec6}, the first reason is to increase the resolution
and make the peak at the reflector's location isotropic. If we do
not sum over equi-distributed illumination angles over the sphere,
we get more of "8-shaped" spot, as shown in Figure~\ref{graphI}.
The second reason is that an average over multiple measurements
increases the stability of the imaging functional with respect to
measurement noise. \end{rem}

\begin{rem}
If we could take an image of the medium in the absence of
reflector before taking the real image, we would be in possession
of the boundary data $\{ u_s - u_s^{(\mu)}, \ x\in \partial
\Omega\}$, and thus we would be able to detect the reflector in a
very noisy background. But in some practical situations
\cite{hsieh2011imaging}, it is not possible to get an image
without the reflector. As it will be shown in section \ref{sec6},
second-harmonic generation can be seen as a powerful contrast
imaging approach \cite{hsieh2011imaging}. In fact, we will prove
that the second harmonic image is much more stable with respect to
the medium noise and to the volume of the particle than the
fundamental frequency image.
\end{rem}

\subsection{Second-harmonic backpropagation}
 If we write a similar imaging functional for the second-harmonic field $v$, assuming that we are in possession
 of the boundary data $\{ v(x), \ x\in \partial\Omega\}$, we
get \begin{equation}\label{DefJ}
 \forall z^S \in \Omega,\ J_\theta (z^S)=\int_{\partial\Omega \times \mathbb{S}^1} v(x)\overline{G_{2\omega}^{(0)}(x,z^S)}
 e^{-2i \omega \theta \cdot z^S} d\sigma(x)d\sigma(\theta).
\end{equation} As before, using (\ref{DAv}) we can expand $J$ in terms of $\delta$ and $\mu$. Considering
first-order terms in $\delta$ and $\mu$ we get
\begin{multline}
 J(z^S)=- \pi \delta^2 \int_{\mathbb{S}^1} e^{-2i\omega \theta \cdot z^S}
 \bigg[ (S)_{det}^\theta \Big(\int_{\partial \Omega} \overline{\G{x}{z^S}}
  \G{x}{z_r} d\sigma(x) \\ - 4 \omega^2 \int_{\partial \Omega} \overline{\G{x}{z^S}}\int_\Omega
 \mu(y) \G{y}{x}  \G{y}{z_r} dy d\sigma(x) \Big) \\  +(S)_{rand}^\theta \int_{\partial \Omega}
  \overline{\G{x}{z^S}} \G{x}{z_r} d\sigma(x)   \bigg] d\sigma(\theta)+
  O(\delta^3 + \delta^2
  ||\mu||_{\mathcal{C}^{0,\alpha^\prime}}^2),
\end{multline}
where $0<\alpha^\prime\leq \alpha$. Now, if we define
$Q_{2\omega}$ as
\begin{equation}\label{DefQ}
 Q_{2\omega}(x,z)= \int_{\partial \Omega} \G{y}{x}  \overline{\G{y}{z}}
 d\sigma(y).
\end{equation}
We have
\begin{multline}\label{DAJ}
 J(z^S)= - \pi \delta^2 \int_{\mathbb{S}^1} e^{-2i\omega \theta \cdot z^S} \bigg[
  (S)_{det}^\theta \left( Q_{2\omega}(z_r,z^S) -
 4\omega^2 \int_{\Omega_\mu}  \mu(y) \G{y}{z_r}  Q_{2\omega}(y,z^S) dy \right) \\+ (S)_{rand}^\theta Q_{2\omega}(z_r,z^S)
  \bigg] d\sigma(\theta) +
 O(\delta^3 + \delta^2 ||\mu||_{\mathcal{C}^{0,\alpha^\prime}}^2 ).
\end{multline}

\section{Statistical analysis} \label{sec6}
In this section, we perform a resolution and stability analysis of
both functionals. Since the image we get is a superposition of a
deterministic image and of a random field created by the medium
noise, we can compute the expectation and the covariance functions
of those fields in order to estimate the signal-to-noise ratio.
For the reader's convenience we give our main results in the
following proposition.

\begin{proposition} Let $l_\mu$ and $\sigma_\mu$  be respectively the
correlation length and the standard deviation of the process
$\mu$. Assume that  $l_\mu$ is smaller than the wavelength
$2\pi/\omega$. Let $(SNR)_I$ and $(SNR)_J$ be defined by
\begin{equation}
\label{snridef}
(SNR)_I=\frac{\mathbb{E}[I(z_r)]}{(Var[I(z_r)])^{1/2}},
\end{equation}
and
\begin{equation}
\label{snrjdef} (SNR)_J= \frac{\mathbb{E}[J(z_r)]
}{(Var[J(z_r)])^{\frac{1}{2}}}.
\end{equation}
We have
\begin{equation}
(SNR)_I \approx \frac{\sqrt{2}\pi^{3/2}\omega\delta^2
U_I}{\sigma_\mu l_\mu \sqrt{\omega\text{ diam } \Omega_\mu}}
\frac{|\sigma_r-1|}{\sigma_r+1},
\end{equation}
and
\begin{equation}
(SNR)_J \geq \frac{  l_\mu^{\alpha} \left(\int_{\mathbb{S}^1}\left(\sum_{k,l} \chi_{k,l} \theta_k
 \theta_l \right) d\theta\right)   }{\sqrt{C} \sigma_\mu \min(\omega^{-\alpha},1)  \max_{k,l} \left\vert \chi_{k,l} \right\vert\sqrt{ \left(\omega \text{diam } \Omega_\mu\right)^{3+2\alpha} + 1 }  }.
\end{equation}
Here, $\text{diam}$ denotes the diameter, $\alpha $ is the upper bound for Holder-regularity of the random
process $\mu$ (see section \ref{noiseassumptions}).
\end{proposition}

\subsection{Assumptions on the random process $\mu$}\label{noiseassumptions}
Let $z(x)$, $x\in \mathbb{R}^2$ be a stationary random process with
Gaussian statistics, zero mean, and a covariance function given by
$R(\vert x-y \vert )$ satisfying $R(0)= \sigma_\mu^2$, $\vert R(0)- R(s)\vert \ \leq \sigma_\mu^2 \frac{s^{2\alpha}}{l_\mu^{2\alpha}} $ and $R$ is decreasing. Then, $z$ is a
$\mathcal{C}^{0,\alpha'}$ process for any $\alpha' <\alpha$
(\cite[Theorem 8.3.2]{adler2010geometry}).  Let $F$ be a smooth
odd bounded function, with derivative bounded by one. For example
 $F=\arctan$ is a suitable choice. Take $$\mu(x)
=F[z(x)].$$ Then $\mu$ is a bounded
$\mathcal{C}^{0,\alpha'}$  stationary process with zero
mean. We want to compute the expectation of its norm. Introduce
 \begin{equation}
p(h)=\max_{\Vert x-y \Vert \leq \sqrt{2} h }  \mathbb{E} \vert z(x) - z(y) \vert.
\end{equation}
One can also write $p(u)=\sqrt{2} \sqrt{R(0) - R(\sqrt{2}u)}$.
According to \cite{adler2010geometry}, for all $h, t \in
\Omega_\mu$, almost surely,
 \begin{equation}
\vert z(t+h) - z(t) \vert \leq  16\sqrt{2 }[\log(B)]^{1/2} p(\frac{\vert h \vert}{l_\mu} ) + 32\sqrt{2} \int_0^{\frac{\vert h \vert}{l_\mu}} \left(-\log u \right)^{1/2} dp(u),
\end{equation} where $B$ is a positive random variable with $\mathbb{E} [B^n] \leq (4\sqrt{2})^n$
(\cite[Formula 3.3.23]{adler2010geometry}). We have that
\begin{equation}
p( \vert h \vert  ) \leq  \sqrt{2}^{1+\alpha} \sigma_\mu  \frac{\vert h \vert^\alpha}{l_\mu^\alpha} .
\end{equation} By integration by parts we find that
\begin{equation}
\int_0^{\frac{\vert h \vert}{l_\mu}} \left(-\log u \right)^{1/2} dp(u) = \left[ (-\log u)^{1/2} p(u)
\right]_0^{\frac{\vert h \vert}{l_\mu}} + \frac{1}{2} \int_0^{\frac{\vert h \vert}{l_\mu} } (-\log u)^{-1/2} u^{-1} p(u) du.
\end{equation}
For any $\varepsilon >0$,  since $s^\varepsilon\sqrt{-\log s} \leq
\frac{1}{\sqrt{ \varepsilon}} e^{1/2} $ on $[0,1]$,
we have, as $\vert h \vert $ goes to $0$, that \begin{equation} \left[ (-\log
u)^{1/2} p(u) \right]_0^{\frac{\vert h \vert}{l_\mu}} \leq
e^{\frac{1}{2 }} \frac{\sqrt{2}^{1+\alpha}\sigma_\mu}{\sqrt{ \varepsilon}}
 \frac{\vert h \vert^{\alpha-\varepsilon}}{l_\mu^\alpha} .
\end{equation}
Similarly, when $\vert h \vert <\frac{1}{2e}$, for every $0<u<\vert h \vert $, $$(-\log
u)^{-1/2} s^{-1} p(u) \leq  \sqrt{2}^{1+\alpha} \sigma_\mu
\frac{u^{\alpha-1}}{l_\mu^\alpha }.$$ So we get, when $\vert h \vert $ goes to $0$, for every
$\varepsilon >0$,
\begin{equation}
\int_0^{\frac{\vert h \vert}{l_\mu}} \left(-\log u \right)^{1/2} dp(u) \leq \frac{e^{\frac{1}{2}}  \sqrt{2}^{1+\alpha} \sigma_\mu }{\sqrt{ \varepsilon}} \frac{ \vert h \vert^{\alpha-\varepsilon}}{l_\mu^\alpha} + \frac{ \sqrt{2}^{1+\alpha} \sigma_\mu }{\alpha }  \frac{\vert h \vert^\alpha}{l_\mu^\alpha}.
\end{equation}
Therefore, when $\vert h \vert $ goes to zero, we have for any $\varepsilon>0$:
\begin{equation}
\vert z(t+h) -z(t) \vert \leq 32 \sqrt{2}^\alpha \log(B)^{1/2} \sigma_\mu \frac{\vert h \vert^\alpha}{l_\mu^\alpha}  +64 e^{\frac{1}{2}} \sqrt{2}^\alpha \sigma_\mu \frac{1}{l_\mu^\alpha} \left[ \frac{1}{\sqrt{\varepsilon}}  \vert h \vert^{\alpha-\varepsilon} +\frac{1}{2} \vert h \vert^\alpha \right].
\end{equation}
Since $F'\leq 1$, composing by $F$ yields, for any $x,y \in
\mathbb{R}^2$,
\begin{equation}
\vert \mu(x) -\mu(y) \vert \leq   \vert z(x)-z(y) \vert .
\end{equation} We get the following estimate on $\Vert \mu \Vert_{\mathcal{C}^{0,\alpha'}}$, for any $\alpha' \in ]0,\alpha[$,
almost surely
\begin{equation}
\sup_{\substack{x, y \in \Omega_\mu \\ \vert x-y \vert \leq h }}\frac{\vert \mu (x) - \mu (y) \vert }{\vert x - y \vert^{\alpha'}} \leq 32 \sqrt{2}^\alpha \log(B)^{1/2} \sigma_\mu \frac{h^{\alpha-\alpha'}}{l_\mu^\alpha} +64 e^{\frac{1}{2}} \sqrt{2}^\alpha \sigma_\mu \frac{1}{l_\mu^\alpha} \left[ \frac{1}{\sqrt{\alpha-\alpha'}}    +\frac{1}{2}  h^{\alpha-\alpha'} \right]
\end{equation}
\begin{equation}
\Vert \mu \Vert_{\mathcal{C}^{0,\alpha'}} \leq 64\sqrt{2}^\alpha \frac{e^{\frac{1}{2}} \left[ \log (B)^{1/2} +1 \right] }{\sqrt{\alpha- \alpha'}}
 \frac{\sigma_\mu}{l_\mu^{\alpha}},
\end{equation} which gives, since $\mathbb{E}[\log B] \leq  \mathbb{E}[B] -1 \leq 4\sqrt{2} -1$
\begin{equation} \label{estc0alpha}
\mathbb{E}[\Vert \mu \Vert_{\mathcal{C}^{0,\alpha'}}^2] \leq  64^2 2^{4+\alpha}  \frac{e}{\alpha- \alpha'}  \frac{\sigma_\mu^2}{l_\mu^{2\alpha}}.
\end{equation}
\subsection{Standard backpropagation}

\subsubsection{Expectation}
We use (\ref{DAI}) and the fact that $\mathbb{E}(\mu)(x)=0, \
\forall x\in \Omega$, to find that
\begin{equation}\label{ExpectI}
\mathbb{E}[I(z^S)]=-2\pi \delta^2\frac{\sigma_r-1}{\sigma_r+1 }  U_I\int_{\mathbb{S}^1}e^{-i\omega \theta \cdot (z^S-z_r) }\theta^\top R_\omega(z^S,z_r) \theta d\theta.
\end{equation}
We now use the Helmholtz-Kirchoff theorem. Since (see \cite{ammarimethods}):
\begin{equation}
\lim_{R\rightarrow \infty} \int_{\vert x\vert
=R}\nabla G_\omega^{(0)}(x,y) \overline{\nabla G_\omega^{(0)}(z,y)}^\top dy = \frac{1}{\omega} \nabla_z \nabla_x \text{ Im}\left[ G_\omega^{(0)} (x,z)\right]
\end{equation}
and \begin{equation}\label{IMG}
\text{ Im}\left[ G_\omega^{(0)}(x,z)\right]=\frac{1}{4}J_0(\omega \vert x-z\vert),
\end{equation} we can compute an approximation of $R_\omega$.
\begin{multline}\label{DAR1}
\frac{1}{\omega} \nabla_z \nabla_x \text{ Im}\left[ G_\omega^{(0)} (x,z)\right] = \frac{1}{4}
\bigg[ \omega J_0(\omega \vert x-z \vert) \left(\frac{(x-z)}{\vert x-z\vert }\frac{(x-z) ^\top}{\vert x-z\vert }
\right) \\ -\frac{2J_1(\omega \vert x-z \vert)}{\vert x-z \vert }\left(\frac{(x-z)}{\vert x-z\vert }
\frac{(x-z) ^\top}{\vert x-z\vert } \right) \\ + \frac{J_1(\omega \vert x-z\vert)}{\vert x-z\vert }
 I_2  \bigg],
\end{multline}
where $I_2$ is the $2\times2$ identity matrix. We can see that
$R_\omega$ decreases as $\vert z_r-z^S\vert^{-\frac{1}{2}}$.  The
imaging functional has a peak at location $z^S=z_r$. Evaluating
$R_\omega$ at $z^S=z_r$ we get
\begin{equation}\
R_\omega (z_r,z_r)= \frac{\omega}{8} I_2.
\end{equation}So we get the expectation of $I$ at point $z_r$:
\begin{equation}\label{ExpectIb}
\mathbb{E}[I(z_r)]\approx  -\frac{\pi^2(\sigma_r-1)}{2(\sigma_r+1)}\omega\delta^2 U_I.
\end{equation}

\subsubsection{Covariance}
Let
\begin{equation} \label{eq41}
\text{Cov}\left(I(z^S),I(z^{S'}) \right) =
\mathbb{E}\bigg[\left(I(z^S)- \mathbb{E}[I(z^S)] \right)
\overline{\left( I(z^{S'})- \mathbb{E}[I(z^{S'})]\right)} \bigg].
\end{equation}
Define \begin{equation} \widetilde{R}_\omega (z^S,z_r,y) =
\int_{\partial \Omega} \overline{\nabla G_\omega^{(0)}(x,z^S)}
\left(\nabla\nabla G_\omega^{(0)}(x,y) \nabla
G_\omega^{(0)}(y,z_r)\right)^\top d\sigma(x).
\end{equation}
Using (\ref{DAI}) and (\ref{ExpectIb}), we get
\begin{multline}
I(z^S)- \mathbb{E}[I(z^S)]= \int_{\partial \Omega \times
\mathbb{S}^1} \frac{1}{i\omega} e^{-i\omega \theta \cdot z^S}\theta^\top \overline{\nabla G_\omega^{(0)}
 (x,z^S)^\top } u_s^{(\mu)}(x) dxd\theta\\ -2\pi \delta^2\frac{\sigma_r-1}{\sigma_r+1 }
  U_I \int_{\mathbb{S}^1} e^{-i\omega \theta \cdot (z^S-z_r) }
   \bigg[\int_\Omega \mu(y) \theta^\top \widetilde{R}_\omega (z^S,z_r,y)\theta dy \bigg]
   d\theta.
\end{multline}
The computations are a bit tedious. For brevity, we write the
quantity above as
\begin{equation}
I(z^S)-\mathbb{E}[I(z^S)] = A_I(z^S) + B_I(z^S),
\end{equation}
with
\begin{equation}
A_I(z^S)= \int_{\partial \Omega\times \mathbb{S}^1}
\frac{1}{i\omega} e^{-i\omega \theta \cdot z^S}\theta^\top
\overline{\nabla G_\omega^{(0)} (x,z^S) } u_s^{(\mu)}(x)
dxd\theta,
\end{equation}
and
\begin{equation}
B_I(z^S) = -2\pi \delta^2\frac{\sigma_r-1}{\sigma_r+1 }  U_I
\int_{\mathbb{S}^1} e^{-i\omega \theta \cdot (z^S-z_r) }
\bigg[\int_\Omega \mu(y) \theta^\top \widetilde{R}_\omega
(z^S,z_r,y)\theta dy \bigg] d\theta.
\end{equation}

We now compute each term of the product in (\ref{eq41})
separately.
\paragraph{Main speckle term:}
We need to estimate the typical size of $A_I$. From
(\ref{estimateborn}), keeping only terms of first-order in $\mu$
yields
\begin{equation}
A_I(z^S)= -\int_{\partial \Omega\times \mathbb{S}^1}
\frac{1}{i\omega} e^{-i\omega \theta \cdot z^S}\theta^\top
\overline{\nabla G_\omega^{(0)} (x,z^S) } \int_\Omega  \mu(y)
\nabla G_\omega^{(0)}(x,y)\cdot \nabla U_0(y)dydxd\theta
+O(\Vert \mu \Vert_\infty^2),
\end{equation}
so we have:
\begin{equation}
A_I(z^S)= -U_I \int_{\Omega\times \mathbb{S}^1}
e^{-i\omega\theta\cdot (z^S-y)} \mu(y) \theta^\top
R_\omega(z^S,y)\theta dyd\theta,
\end{equation}
and hence,
\begin{multline}
A_I (z^S)\overline{A_I (z^{S'})} =U_I^2 \int_{\mathbb{S}^1}
e^{-i\omega \theta \cdot (z^S-z^{S'})} \\  \bigg[ \int
\int_{\Omega \times \Omega} e^{i\omega \theta \cdot (y-y')} \mu(y)
\mu(y') \theta^\top R_\omega(z^S,y) \overline{R_\omega(z^{S'},y')}
\theta dydy'\bigg]d \theta.
\end{multline}
We assume that the medium noise is localized and stationary on its
support $\Omega_\mu$. We also assume that the correlation length
$l_\mu$ is smaller than the wavelength. We  note $\sigma_\mu$ the
standard deviation of the process $\mu$. We can then write:
\begin{equation}
\mathbb{E}\bigg[A_I (z^S)\overline{A_I (z^{S'})}\bigg] = U_I^2
\sigma_\mu^2 l_\mu^2\int_{\mathbb{S}^1}e^{i\omega \theta \cdot
(z^S-z^{S'}) } \int_{\Omega _\mu} \theta^\top R_\omega (z^S,y)
\overline{R_\omega (z^{S'},y)}\theta dy d\theta.
\end{equation}
We introduce
\begin{equation}\label{defP}
P_\omega(z^S,y,z^{S'}):=\int_{\mathbb{S}^1}e^{i\omega \theta \cdot
(z^S-z^{S'}) } \theta^\top R_\omega (z^S,y)   \overline{R_\omega
(z^{S'},y)}\theta d\theta,
\end{equation}
where $R_\omega$ is defined by (\ref{DEFR}). Therefore, we have
\begin{equation}\label{EstA}
\mathbb{E}\bigg[A_I (z^S)\overline{A_I (z^{S'})}\bigg] = U_I^2
\sigma_\mu^2 l_\mu^2 \int_{\Omega _\mu} P_\omega (z^S,y,z^{S'})
dy.
\end{equation}
Hence, $A_I$ is a complex field with Gaussian statistics of mean
zero and covariance given by  (\ref{EstA}). It is a speckle field
and is not localized.

We compute its typical size at point $z^S=z^{S'}=z_r$, in order to
get signal-to-noise estimates. Using (\ref{DAR1}), we get that for $\vert x-z\vert >>1$:
\begin{equation*}
\lim_{R\rightarrow \infty} \int_{\vert x\vert =R}\nabla
G_\omega^{(0)}(x,y) \overline{\nabla G_\omega^{(0)}(z,y)}^\top dy
= \frac{\omega }{4}J_0(\omega \vert x-z\vert)
\left(\frac{(x-z)}{\vert x-z\vert }\frac{(x-z) ^\top}{\vert
x-z\vert } \right).
\end{equation*}
Since we have, for $\vert x-z\vert >> 1$, \begin{equation}\label{equivJ0}
 J_0(\omega \vert x -z \vert ) \sim \frac{\sqrt{2} \cos(\omega \vert x-z \vert -
 \frac{\pi}{4})}{\sqrt{\pi \omega \vert x-z\vert}},
\end{equation} we obtain that \begin{equation}\label{DAR}
R_\omega (x,z)  \approx \frac{\sqrt{\omega}\cos (\omega \vert x-z\vert -\pi/4)}{2\sqrt{2\pi }} \vert x-z\vert ^{-1/2}   \left(\frac{(x-z)}{\vert x-z\vert }\frac{(x-z) ^\top}{\vert x-z\vert } \right) \ \text{for } \vert x-z\vert >>1.
\end{equation}
Now we can write
\begin{equation}
\mathbb{E}\bigg[A_I (z_r)\overline{A_I (z_r)}\bigg]  \approx U_I^2
\sigma_\mu^2 l_\mu^2 \int_{\Omega _\mu}
\left(\frac{\sqrt{\omega}}{2\sqrt{2\pi }}\right)^2\frac{1}{2} \vert
y-z_r\vert ^{-1} \int_{\mathbb{S}^1} \theta^\top
\left(\frac{(y-z_r)}{\vert y-z_r\vert }\frac{(y-z_r) ^\top}{\vert
y-z_r\vert }\right)\theta d\theta  dy.
\end{equation}
If we compute the term:
\begin{equation}
\int_{\mathbb{S}^1} \theta^\top \left(\frac{(y-z_r)}{\vert y-z_r\vert }\frac{(y-z_r) ^\top}{\vert y-z_r\vert }\right)\theta d\theta = \int_0^{2\pi} \bigg[\left( \frac{(y-z_r)_1}{\vert y-z_r \vert}\right)^2 \cos^2\theta + \left(\frac{(y-z_r)_2}{\vert y-z_r \vert}\right)^2 \sin^2\theta\bigg] d\theta,
\end{equation}
then, after linearization and integration, we get
\begin{equation}\label{eqpi}
\int_{\mathbb{S}^1} \theta^\top \left(\frac{(y-z_r)}{\vert
y-z_r\vert }\frac{(y-z_r) ^\top}{\vert y-z_r\vert }\right)\theta
d\theta = \pi.
\end{equation}
So we have:
\begin{equation}
\mathbb{E}\bigg[A_I (z_r)\overline{A_I (z_r)}\bigg]  \approx  \pi
U_I^2 \sigma_\mu^2 l_\mu^2 \int_{\Omega _\mu}
\left(\frac{\sqrt{\omega}}{4\sqrt{\pi }}\right)^2 \vert
y-z_r\vert ^{-1}  dy,
\end{equation}
and therefore,
\begin{equation}\label{TERMA}
\mathbb{E}\bigg[A_I (z_r)\overline{A_I (z_r)}\bigg]  \approx \pi \frac{\omega }{8} U_I^2
\sigma_\mu^2 l_\mu^2 \text{diam } \Omega_\mu.
\end{equation}
\paragraph{Secondary speckle term:}

We have
\begin{multline}
B_I (z^S)\overline{B_I (z^{S'})}=\left(2\pi \delta^2\frac{\sigma_r-1}{\sigma_r+1 }
U_I \right)^2 \int_{\mathbb{S}^1} e^{-i\omega \theta \cdot (z^S-z^{S'}) }
\\\bigg[\int_\Omega \mu(y) \mu(y') \theta^\top \widetilde{R}_\omega (z^S,z_r,y)
\overline{\widetilde{R}_\omega (z^{S'},z_r,y')}\theta dy dy'\bigg]
d\theta.
\end{multline}
So we get the expectation:
\begin{multline}
\label{DAspeckle2} \mathbb{E}\bigg[B_I (z^S)\overline{B_I
(z^{S'})}\bigg] =\left(2\pi \delta^2\frac{\sigma_r-1}{\sigma_r+1 }  U_I \right)^2 \sigma_\mu^2
l_\mu^2 \\ \int_{\mathbb{S}^1} e^{-i\omega \theta \cdot (z^S-z^{S'}) } \theta^\top\bigg[ \int_{\Omega_\mu}  \widetilde{R}_\omega(z^S,z_r,y)\overline{\widetilde{R}_\omega(z^{S'},z_r,y)} dy \bigg]\theta d\theta.
\end{multline} This term also creates a speckle field on the image.
As before, we compute the typical size of this term at point $z_r$.
We first get an estimate on $\widetilde{R}_\omega$.
\begin{equation}
\vert\left( \widetilde{R}_\omega(z^S,z_r,y)\right)_{i,j} \vert \leq  \vert \partial_j G_\omega^{(0)}(y,z_r) \vert  \vert \sum_{k=1,2}\int_{\partial \Omega} \partial_{y_i}\overline{G_\omega^{(0)}(x,z^S)} \partial_{y_i} \partial_{y_k} G_\omega^{(0)}(x,y) d\sigma(x) \vert.
\end{equation}
We recall the Helmholtz-Kirchoff theorem
\begin{equation}
\int_{\partial \Omega} \overline{G_\omega^{(0)}(x,y)}
G_\omega^{(0)}(x,z) d\sigma(x) \sim \frac{1}{4\omega} J_0(\omega
\vert y-z\vert) \quad \mbox{as } R\rightarrow\infty,
\end{equation}
from which
\begin{equation}
\int_{\partial \Omega} \partial_{y_i}\overline{G_\omega^{(0)}(x,z^S)} \partial_{y_i} \partial_{y_k} G_\omega^{(0)}(x,y) d\sigma(x) = \frac{1}{4\omega} \left(\partial_i \partial_i \partial_k f\right) (z^S-y),
\end{equation}
where $f$ is defined by $f(x)=J_0(\omega \vert x \vert)$.
We have
\begin{equation}
\partial_i \partial_j \partial_k f (x) = \omega \left( \frac{3\left(a_{i,j,k}(x)-b_{i,j,k}(x)\right)}{\vert x \vert^2 }\left[ J_0'(\omega \vert x \vert) -\omega \vert x \vert  J_0''(\omega \vert x \vert) \right]+ a_{i,j,k}(x) \omega^2 J_0^{(3)}(\omega \vert x \vert )\right),
\end{equation} where $a_{i,j,k}$ and $b_{i,j,k}$ are rational fractions in the coefficients of $x$ bounded by $1$.
Now, recall the power series of $J_0$:
\begin{equation}
J_0(z) = \sum_k (-1)^k \frac{\left(\frac{1}{4}z^2\right)^k}{(k!)^2}.
\end{equation} We can write
\begin{equation}
 J_0'(\omega \vert x \vert) -\omega \vert x \vert  J_0''(\omega \vert x \vert)
 =-\frac{\omega^3}{4} \vert x\vert^3 +o(\vert x \vert^3).
\end{equation}
Hence, since $J_0^{(3)}(x)\sim \frac{3}{4} x$ when $x\rightarrow
0$, we can prove the following estimate for $x$ around $0$:
\begin{equation}
\frac{1}{4\omega}(\partial_i \partial_j \partial_k f)(x)  \sim  \frac{3b_{i,j,k}(x)}{16} \omega^3 \vert x \vert .
\end{equation}
In order to get the decay of $\widetilde{R}_\omega$ for large
arguments we use the following formulas: $J_0'=-J_1$,
$J_0''=\frac{1}{x}J_1-J_0$, and $J_0^{(3)}= J_1-\frac{1}{x^2}J_1+
\frac{1}{x} J_0$. We get
\begin{equation}
\frac{1}{4\omega}\vert \partial_i \partial_j \partial_k f (x)
\vert \leq \omega^2 (\omega \vert x \vert )^{-1/2} \quad \mbox{as
} x\rightarrow \infty.
\end{equation}
We also have the following estimate:
\begin{equation}
\vert \nabla G_\omega^{(0)}(y,z_r)
\vert \leq \left(\frac{2}{\pi}\right)^{1/2}\max \left( \frac{1}{ \vert y-z_r \vert},
 \frac{\omega}{\sqrt{\omega \vert y-z_r \vert }}\right).
\end{equation}
We can now write the estimate on $\widetilde{R_\omega}_{i,j}$
\begin{equation}
\vert \widetilde{R_\Omega}(z^S,z_r,y)_{i,j}\vert  \leq \omega^2 \left(\frac{2}{\pi}\right)^{1/2} \min\left(\omega \vert y-z_r \vert , \frac{1}{\sqrt{\omega \vert y-z^S \vert}} \right) \max \left( \frac{1}{\omega \vert y-z_r \vert}, \frac{1}{\sqrt{\omega \vert y-z_r \vert }}\right).
\end{equation}We can now go back to estimating the term $B_I$. We split the domain of integration $\Omega_\mu =B(z_r,\omega^{-1})
\cup \Omega_\mu \backslash B(z_r,\omega^{-1})$ to get
\begin{multline}
\left\vert \mathbb{E}\bigg[B_I(z_r)\overline{B_I(z_r)}\bigg] \right\vert  \leq    \left(2\pi \delta^2\frac{\sigma_r-1}{\sigma_r+1 }  U_I \right)^2 \sigma_\mu^2
l_\mu^2 \\ 4\pi \omega^4 \frac{2}{\pi}\bigg[ \int_{\Omega_\mu \backslash B(z_r,\omega^{-1})} \frac{1}{\vert y-z_r\vert^2}dy + \int_{B(z_r,\omega^{-1})} \omega^2 figuresdy \bigg].
\end{multline}
Hence,
\begin{equation}\label{TERMB}
\left\vert \mathbb{E}\bigg[B_I(z_r)\overline{B_I(z_r)}\bigg] \right\vert \leq 8\left(2\pi \delta^2\frac{\sigma_r-1}{\sigma_r+1 }  U_I \right)^2  \omega^4 \sigma_\mu^2
l_\mu^2 \log(\omega \text{ diam }\Omega_\mu ).
\end{equation}

\paragraph{Double products:}
The double products $A_I \overline{B_I}$ and $B_I \overline{A_I}$  have a typical amplitude that is the geometric mean of the
typical amplitudes of $A_I$ and $B_I$. So they are always
smaller than one of the main terms $\vert A_I\vert^2$ or $\vert
B_I\vert^2$.
\subsubsection{Signal-to-noise ratio estimates}
We can now give an estimate of the signal-to-noise ratio $(SNR)_I$
defined by (\ref{snridef}). Using (\ref{ExpectIb}), (\ref{TERMA}), and
(\ref{TERMB})  we get
\begin{equation}\label{SNRI1}
(SNR)_I \approx \frac{\frac{\pi^2(\sigma_r-1)}{2(\sigma_r+1)}\omega\delta^2 U_I}{\sigma_\mu l_\mu \left( \pi
\frac{\omega}{8} \text{ diam }\Omega_\mu
+8\left(2\pi \delta^2\frac{\sigma_r-1}{\sigma_r+1 }  U_I \right)^2  \omega^4  \log(\omega \text{ diam }\Omega_\mu )\right)^{1/2}},
\end{equation}
Since $\delta << \frac{2\pi}{\omega}$ we have that $\delta \omega
<< 1$, so we can estimate $(SNR)_I$ as follows
\begin{equation}\label{SNRI}
(SNR)_I \approx \frac{\sqrt{2}\pi^{3/2}\frac{\sigma_r-1}{\sigma_r+1}\omega\delta^2 U_I}{\sigma_\mu l_\mu \sqrt{\omega \text{ diam }
\Omega_\mu}}.
\end{equation}
The perturbation in the image $I$ comes from different phenomena.
The first one, and the most important is the fact that we image
not only the field scattered by the reflector, but also the field
scattered by the medium's random inhomogeneities. This is why the
signal-to-noise ratio depends on the volume and the contrast of
the  particle we are trying to locate. It has to stand out from
the background. The other terms in the estimate (\ref{SNRI1}) of
$(SNR)_I$ are due to the phase perturbation of the field scattered
by the particle when it reaches the boundary of $\Omega$  which
can be seen as a travel time fluctuation of the scattered wave by
the reflector. Both the terms are much smaller than the first one.
$(SNR)_I$ depends on the ratio ${\omega}/{l_\mu}$. If the medium
noise has a shorter correlation length, then the perturbation
induced in the phase of the fields will more likely self average.

\subsection{Second-harmonic backpropagation}
\subsubsection{Expectation}
We have:
\begin{multline}
\mathbb{E}[J(z^S)] = -\pi \delta^2  \int_{\mathbb{S}^1} e^{-2i\omega \theta \cdot z^S} \bigg[
(S)_{det}^\theta\int_{\partial \Omega} \overline{\G{x}{z^S}} \G{x}{z_r} dx \\
+ \mathbb{E}[(S)_{rand}^\theta] \int_{\partial \Omega}
\overline{\G{x}{z^S}} \G{x}{z_r} dx \bigg]d\theta.
\end{multline}
Since $\mathbb{E}[(S)_{rand}^\theta] =0$ we obtain  by using
(\ref{DEFSDET}) that
\begin{equation}
 \mathbb{E}[J(z^S)] =
 \pi \delta^2 \omega^2 U_I^2 \int_{\mathbb{S}^1} \left(\sum_{k,l} \chi_{k,l} \theta_k \theta_l \right) e^{2i\omega \theta \cdot (z_r-z^S)} d\theta \int_{\partial \Omega} \overline{\G{x}{z^S}} \G{x}{z_r}
 dx.
\end{equation}
If we define $\widetilde{Q}_{2\omega}$ as
\begin{equation}\label{defqtilde}
\widetilde{Q}_{2\omega}(x,y)= \int_{\mathbb{S}^1} \left(\sum_{k,l} \chi_{k,l} \theta_k
 \theta_l \right) e^{2i\omega \theta \cdot (x-y)} d\theta,
\end{equation}
then it follows that
\begin{equation}
 \mathbb{E}[J(z^S)] = \delta^2 \omega^2 U_I^2 \widetilde{Q}_{2\omega}(z_r,z^S) Q_{2\omega}(z_r,z^S),
\end{equation} where $Q_{2\omega}$ is given by (\ref{DefQ}).
 To get the typical size of this term we first
use the Helmholtz-Kirchhoff theorem \cite{ammarimethods}:
\begin{equation}\label{EQUIVQ}
Q_{2\omega}(z_r,z^S) \sim \frac{1}{2\omega}\text{ Im} \left( \G{z_r}{z^S}\right).
\end{equation}
Therefore,  we obtain that
\begin{equation}\label{ExpectJ}
 \mathbb{E}[J(z_r)] = \frac{\pi}{8} \delta^2\omega U_I^2  \int_{\mathbb{S}^1} \left(\sum_{k,l} \chi_{k,l} \theta_k
 \theta_l \right) d\theta.
\end{equation}

\subsubsection{Covariance}
We have:
\begin{multline}
J(z^S)-\mathbb{E}[J](z^S) = \pi \delta^2 \int_{\mathbb{S}^1}e^{-2i\omega \theta \cdot z^S}\Big[ (S)_{det}^\theta 4\omega^2 \int_\Omega \G{s}{z_r}\mu(s) Q_{2\omega}(s,z^S)ds \\ - (S)_{rand}^\theta Q_{2\omega}(z_r,z^S) \Big] d\theta.
\end{multline}
Denote by \begin{equation} A_J(z^S)= 4 \pi \delta^2 \omega^2 \int_{\mathbb{S}^1}e^{-2i\omega \theta \cdot z^S} (S)_{det}^\theta \int_\Omega
\G{s}{z_r}\mu(s) Q_{2\omega}(s,z^S)ds d\theta,
\end{equation} and
\begin{equation}B_J(z^S)=\pi \delta^2 \int_{\mathbb{S}^1}e^{-2i\omega \theta \cdot z^S} (S)_{rand}^\theta Q_{2\omega}(z_r,z^S) d\theta.
 \end{equation}
Then we can write the covariance function,
\begin{equation}
\text{Cov}\left(J(z^S), J(z^{S'}) \right) = \mathbb{E}\bigg[\left(J(z^S)- \mathbb{E}[J(z^S)] \right)
\overline{\left( J(z^{S'})- \mathbb{E}[J(z^{S'})]\right)} \bigg],
\end{equation}
in the form
\begin{equation}
\text{Cov}\left(J(z^S), J(z^{S'}) \right) =\mathbb{E}\bigg[ A(z^S)
\overline{A(z^{S'})} +B(z^S)\overline{B_J(z^{S'})}
+A_J(z^S)\overline{B_J(z^{S'})} +
\overline{A_J(z^S)}B_J(z^{S'})\bigg].
\end{equation}
We will now compute the first two terms separately and then we deal with the double products.

\paragraph{The speckle term $A_J\overline{A_J}$:}

\medskip
From
\begin{multline}
A_J(z^S) \overline{A_J(z^{S'})} = 16 \pi^2\delta^4\omega^{4} \int_{\mathbb{S}^1} e^{-2i\omega \theta \cdot( z^S-z^{S'})}\vert (S)_{det}^\theta \vert^2
 \\ \int \int _{\Omega \times \Omega} \G{s}{z_r}
\overline{\G{s'}{z_r}} \mu(s) \overline{\mu(s')}
Q_{2\omega}(s,z^S) \overline{Q_{2\omega}(s',z^{S'})} ds ds' d\theta,
\end{multline}
it follows by using (\ref{DEFSDET}) that
\begin{multline}
A_J(z^S) \overline{A_J(z^{S'})} =16 \pi^2\delta^4 \omega^{8} U_I^4 \int_{\mathbb{S}^1} e^{-2i\omega \theta \cdot( z^S-z^{S'})} \vert \sum_{k,l} \chi_{k,l} \theta_k \theta_l  \vert^2 d\theta \\ \int
\int _{\Omega \times \Omega} \G{s}{z_r} \overline{\G{s'}{z_r}} \mu(s)
 \overline{\mu(s')} Q_{2\omega}(s,z^S) \overline{Q_{2\omega}(s',z^{S'})} ds
 ds'.
\end{multline}
If we write $C_\mu(s,s')=\mathbb{E}[\mu(s) \mu(s')]$, then we find
that
\begin{multline}
\mathbb{E}[A_J(z^S) \overline{A_J(z^{S'})}] =  16 \pi^2\delta^4 \omega^{8} U_I^4  \int_{\mathbb{S}^1}e^{-2i\omega \theta \cdot( z^S-z^{S'})} \vert \sum_{k,l} \chi_{k,l} \theta_k \theta_l  \vert^2 d\theta \\
\int \int _{\Omega \times \Omega} \G{s}{z_r} \overline{\G{s'}{z_r}} C_\mu(s,s') Q_{2\omega}(s,z^S)
\overline{Q_{2\omega}(s',z^{S'})} ds ds' ,
\end{multline}
since $\mu$ is real.

As previously, we assume that the medium noise is localized and
stationary on its support (which is $\Omega_\mu$).  We note
$\sigma_\mu$ the standard deviation of the process $\mu$ and
$l_\mu$ its correlation length. We can then write
\begin{multline}
\mathbb{E}[A_J(z^S) \overline{A_J(z^{S'})}] =
 16 \pi^2\delta^4 \omega^{8} U_I^4 \sigma_\mu^2 l_\mu^2  \int_{\mathbb{S}^1}e^{-2i\omega \theta \cdot( z^S-z^{S'})} \vert \sum_{k,l} \chi_{k,l} \theta_k \theta_l  \vert^2 d\theta \\ \int _{\Omega_\mu} \vert \G{s}{z_r}\vert^2 Q_{2\omega}(s,z^S)
\overline{Q_{2\omega}(s,z^{S'})} ds .
\end{multline}
The term  $\mathbb{E}[A_J(z^S) \overline{A_J(z^{S'})}]$ shows the generation of a non localized
speckle image, creating random secondary peaks. We will later
estimate the size of those peaks in order to find the
signal-to-noise ratio. We compute the typical size of this term.
We get, using (\ref{EQUIVQ}):
\begin{multline} \label{CovA}
\mathbb{E}[A_J(z^S) \overline{A_J(z^{S'})}] \approx 4\pi^2 U_I^4\delta^4 \omega^{6} \sigma_\mu^2 l_\mu^2 \\ \int_{\mathbb{S}^1}\vert \sum_{k,l}
\chi_{k,l} \theta_k \theta_l  \vert^2d\theta  \int
_{\Omega_\mu} \vert \G{s}{z_r}\vert^2 \text{ Im } \G{s}{z^S}
\text{ Im }  \G{s}{z^{S'}} ds .
\end{multline}
Then we use the facts that  $$\vert \G{x}{y} \vert  \approx \frac{1}{4\sqrt{\pi  2\omega}} \vert
x-y \vert^{-1/2}$$ and $$\text{ Im }\G{x}{y} = \frac{1}{4}
J_0(2\omega \vert x-y \vert )\approx \frac{\cos \left( 2\omega \vert x-y \vert - \pi/4 \right) }{4 \sqrt{\pi \omega}}
\vert x-y\vert ^{-1/2}$$ if $\vert x-y\vert >> 1.$ Then, as
previously, we write $\Omega_\mu= \Omega_\mu \backslash
B(z_r,\omega^{-1}) \cup B(z_r,\omega^{-1})$. Using (\ref{CovA}), we arrive at
\begin{multline}
\mathbb{E}[A_J(z_r) \overline{A_J(z_r)}] \approx 4\pi^2 U_I^4\delta^4 \omega^{6} \sigma_\mu^2 l_\mu^2\int_{\mathbb{S}^1} \vert \sum_{k,l} \chi_{k,l} \theta_k \theta_l  \vert^2 d\theta \\   \bigg( \frac{1}{512\pi^2\omega^2}
 \int_{\Omega_\mu \backslash B(z_r,\omega^{-1})}\frac{ \cos^2 \left( 2\omega \vert s-z_r \vert - \pi/4 \right) }{  \vert s-z_r\vert^{2}} ds +\frac{1}{16}\int_{ B(z_r,\omega^{-1})} \vert G_{2\omega}^{(0)}(s,z_r) \vert^2 J_0(2\omega \vert s-z_r \vert )^2ds\bigg),
\end{multline}
which yields
\begin{equation}\label{estAJ}
\mathbb{E}[A_J(z_r) \overline{A_J(z_r)}] \approx \frac{\pi}{128} U_I^4\delta^4\omega^4\sigma_\mu^2 l_\mu^2\log(\omega \text{ diam }\Omega_\mu )\int_{\mathbb{S}^1}  \vert \sum_{k,l} \chi_{k,l} \theta_k \theta_l  \vert^2 d\theta   .
\end{equation}

\paragraph{The localized term $B_J\overline{B_J}$:}

\medskip
We have
\begin{equation}
 B_J(z^S) \overline{B_J(z^{S'})} = \pi^2\delta^4  Q_{2\omega}(z_r,z^S)
 \overline{Q_{2\omega}(z_r, z^{S'})}\int_{\mathbb{S}^1} e^{-2i\omega \theta
 \cdot( z^S-z^{S'})} \vert (S)_{rand}^\theta \vert ^2 d\theta.
\end{equation}
Using (\ref{DEFSRAND})  and (\ref{transfoSrand})  we have that
$(S)^\theta_{rand}$ can be re-written as
\begin{multline}
(S)^\theta_{rand}=- \omega^2 U_I^2 \int_\Omega  \left(\mu(y) e^{i\omega
\theta \cdot y} - \mu(z_r) e^{i\omega \theta \cdot z_r}\right) \\
 \bigg[ \sum_{k,l} \chi_{k,l} \left(\theta_k \theta \cdot \nabla \partial_{x_l}
 G^{(0)}_\omega (z_r,y)+ \theta_l \theta \cdot \nabla \partial_{x_k} G^{(0)}_\omega (z_r,y) \right) \bigg]
 dy.
\end{multline}
We need to get an estimate on $S^\theta_{rand}$'s variance. As in
section \ref{sec2} we have the following estimate for any
$0<\alpha' <1/2$:
\begin{equation}
\frac{1}{4}\vert y-z_r \vert^{\alpha'} \left\vert  \partial_{x_k} \partial_{x_l} H_0^1(\omega \vert y-z_r \vert ) \right\vert \leq \frac{1}{2} \min \left( 1, \sqrt{\frac{2}{\pi}} \omega^{3/2} \vert y-z_r \vert^{\alpha'-1/2} \right) \max\left( 1, \vert y-z_r \vert^{\alpha'-2} \right).
\end{equation}
We get, for any $\alpha' < \min(\alpha, \frac{1}{2})$,
\begin{equation}
\vert S^\theta_{rand} \vert \leq \omega^2 U_I^2\Vert \mu
\Vert_{\mathcal{C}^{0,\alpha'}} \max_{k,l} \left\vert \chi_{k,l}
\right\vert\omega^{2-2\alpha'} \bigg[
\frac{8\sqrt{2\pi}}{3/2+\alpha'} \left( \omega \text{diam }
\Omega_\mu \right)^{3/2+\alpha'} + \frac{\pi}{\alpha'} \bigg],
\end{equation}
and
\begin{multline}
\left\vert  \mathbb{E}[B_J(z^S)\overline{B_J(z^{S'})}] \right\vert \leq \frac{128 \pi^3}{(3/2+\alpha')^2} \omega^{4-2\alpha'} \delta^4   U_I^4\max_{k,l} \left\vert \chi_{k,l} \right\vert ^2 \mathbb{E}\left[ \Vert \mu \Vert_{\mathcal{C}^{0,\alpha'}}^2\right] \\ \bigg[\left( \omega \text{diam } \Omega_\mu \right)^{3+2\alpha'} + \frac{1}{\alpha'} \bigg]  Q_{2\omega}(z_r,z^S)
 \overline{Q_{2\omega}(z_r, z^{S'})}.
\end{multline}
Note that $Q_{2\omega}(z_r,z^S)$, defined in (\ref{DefQ}), behaves
like $\frac{1}{8\omega}J_0(2\omega \vert z_r - z^S \vert)$ which
decreases like $\vert z_r -z^S \vert^{-1/2}$ as $\vert z_r -z^S
\vert$ becomes large. The term $B_J$ is localized around $z_r$. It
may shift, lower or blur the main peak but it will not contribute
to the speckle field on the image. We still need to estimate its
typical size at point $z_r$ in order to get the signal-to-noise
ratio at point $z_r$. Using (\ref{EQUIVQ}) and (\ref{estc0alpha})
we get
\begin{equation}
\mathbb{E}[B_J(z_r)\overline{B_J(z_r)}]    \leq \frac{2^{17+\alpha}   \pi^3}{(3/2+\alpha')^2} \frac{e}{\alpha-\alpha'} \omega^{2-2\alpha'} \delta^4 U_I^4  \max_{k,l} \left\vert \chi_{k,l} \right\vert ^2 \bigg[\left( \omega \text{diam } \Omega_\mu \right)^{3+2\alpha'} + \frac{1}{\alpha'} \bigg] \frac{\sigma_\mu^2}{l_\mu^{2\alpha}}.
\end{equation}
We can write $(\omega \text{diam } \Omega_\mu )^{3+2\alpha' } \leq (\omega \text{diam } \Omega_\mu)^{3+2\alpha} +1$. We can take  $\alpha'=\frac{\alpha}{2}$. Let $C = \frac{2^{18+1/2} \pi^3 e }{(3/2)^2}$.
We get that 
\begin{equation}\label{estBJ}
\mathbb{E}[B_J(z_r)\overline{B_J(z_r)}]    \leq C \omega^2\min \left(\omega^{-2\alpha},1\right)  \delta^4 U_I^4  \max_{k,l} \left\vert \chi_{k,l} \right\vert ^2 \frac{\sigma_\mu^2}{l_\mu^{2\alpha}} \bigg[ \left( \omega \text{diam } \Omega_\mu \right)^{3+2\alpha} + 1 \bigg]   .
\end{equation}
\begin{rem}
We note that even though the term $B_J$ is localized, meaning it
would not create too much of a speckle far away from the
reflector, it is still the dominant term of the speckle field
around the reflector's location.
\end{rem}

\paragraph{The double products $A_J\overline{B_J}$ and $\overline{A_J}B_J$:}
\medskip

This third term has the size of the geometric mean of the first
two terms $A_J$ and $B_J$. So we only need to
concentrate on the first two terms. Also this term is still
localized because of  $Q(z_r,z^S)$ that decreases as $\vert z_r -
z^S\vert^{-1/2}$.

\subsubsection{Signal-to-noise ratio}

 As before,  we define the signal-to-noise ratio $(SNR)_J$ by
 (\ref{snrjdef}).
Using (\ref{ExpectJ}), (\ref{estAJ}) and (\ref{estBJ}),
\begin{equation}
\frac{\mathbb{E}[J(z_r)]}{(Var(J(z_r))^{\frac{1}{2}}}\geq \frac{  l_\mu^{\alpha} \left(\int_{\mathbb{S}^1}\left(\sum_{k,l} \chi_{k,l} \theta_k
 \theta_l \right) d\theta\right)   }{\sqrt{C} \sigma_\mu \min(\omega^{-\alpha},1)  \max_{k,l} \left\vert \chi_{k,l} \right\vert\sqrt{ \left(\omega \text{diam } \Omega_\mu\right)^{3+2\alpha} + 1 }  }.
\end{equation}
The difference here with the standard backpropagation is that the
$(SNR)$ does not depend on neither the dielectric contrast of the
particle, the nonlinear susceptibility nor even the particle's
volume. All the background noise created by the propagation of the
illuminating wave in the medium is filtered because the small
inhomogeneities only scatter waves at frequency $\omega$. The
nanoparticle is the only source at frequency $2\omega$ so it does
not need to stand out from the background. The perturbations seen
on the image $J$ are due to travel time fluctuations of the wave
scattered by the nanoparticle (for the speckle field) and to the
perturbations of the source field at the localization of the
reflector (for the localized perturbation). The second-harmonic
image is more resolved than the fundamental frequency image.

\subsection{Stability with respect to measurement noise}
We now compute the signal-to-noise ratio in the presence of
measurement noise without any medium noise ($\mu=0$). The signal
$u_s$ and $v$ are corrupted by an additive noise $\nu(x)$ on
$\partial \Omega$. In real situations it is of course impossible
to achieve measurements for an infinity of plane waves
illuminations.  So in this part we assume that the functional $J$
is calculated as an average over $n$ different illuminations,
uniformly distributed in $\mathbb{S}^1$. We consider, for each
$j\in [0,n]$, an independent and identically distributed random
process $\nu^{(j)}(x),\ x\in
\partial \Omega$ representing the measurement noise. We use the
model of \cite{TDerivativ}: if we assume that the surface of
$\Omega$ is covered with sensors half a wavelength apart and that
the additive noise has variance $\sigma$ and is independent from
one sensor to another one, we can model the additive noise process
by a Gaussian white noise with covariance function:
$$\mathbb{E}(\nu(x) \overline{\nu(x')}) = \sigma_\nu^2
\delta(x-x'),$$ where $\sigma_\nu = \sigma^2 \frac{\lambda}{2}$.
\subsubsection{Standard backpropagation}
We write, for each $j\in [0,n]$, $u_s^{(j)}$ as
\begin{equation}
u_s^{(j)}(x)=-2\pi \delta^2 \frac{\sigma_r-1}{\sigma_r+1}  U_I e^{i \omega \theta^{(j)} \cdot z_r}
 \nabla G_\omega^{(0)}(x,z_r) \cdot(i\omega \theta^{(j)} )+ o(\delta^2)+ \nu^{(j)}(x),
\end{equation} where $\nu^{(j)}$ is the measurement noise associated with the $j$-th illumination.
We can write $I$ as
\begin{equation}
I(z^S) =\frac{1}{n} \sum_{j=1}^n\int_{\partial \Omega}\frac{1}{i\omega}e^{-i \omega \theta^{(j)} \cdot z^S}(\theta^{(j)})^\top\overline{\nabla G_\omega^{(0)}(x,z^S)}  u_s(x)dx,
\end{equation}
Further, \begin{multline} I(z^S)= -2\pi \delta^2 \frac{\sigma_r-1}{\sigma_r+1}  U_I \frac{1}{n}\sum_{j=1}^n e^{i \omega \theta^{(j)}
\cdot (z_r-z^S)}(\theta^{(j)})^\top R_\omega(z_r,z^S) \theta^{(j)}
\\ + \frac{1}{n}\sum_{j=1}^n \int_{\partial \Omega}\frac{1}{i\omega} e^{-i \omega
\theta^{(j)} \cdot z^S}(\theta^{(j)})^\top\overline{\nabla
G_\omega^{(0)}(x,z^S)} \nu^{(j)} (x)dx.
\end{multline}
We get that
\begin{equation}
\mathbb{E}[I(z^S)]= -2\pi \delta^2 \frac{\sigma_r-1}{\sigma_r+1}  U_I \frac{1}{n}\sum_{j=1}^n e^{i \omega \theta^{(j)}
\cdot (z_r-z^S)}(\theta^{(j)})^\top R_\omega(z_r,z^S) \theta^{(j)},
\end{equation}
so that, using (\ref{DAR1}) and (\ref{IMG})
\begin{equation}
\mathbb{E}[I(z_r)]\sim -\frac{\pi(\sigma_r-1)}{4(\sigma_r+1)}\omega\delta^2 U_I.
\end{equation}
We compute the covariance
\begin{multline}
Cov(I(z^S),I(z^{S'})) = \mathbb{E}\bigg[\frac{1}{n^2}\left( \sum_{j=1}^n \frac{1}{i\omega}e^{-i\omega \theta^{(j)} \cdot
z^S} \int_{\partial\Omega } \nu^{(j)}(x) (\theta^{(j)})^\top
\overline{\nabla G_\omega^{(0)}(x,z^S)} dx\right) \\ \left(\sum_{l=1}^n\frac{-1}{i\omega} e^{i\omega \theta^{(l)} \cdot
z^{S'}} \int_{\partial\Omega} \nu^{(l)}(x') (\theta^{(l)})^\top\nabla G_\omega^{(0)}(x',z^{S'})dx' \right) \bigg],
\end{multline}
and obtain that
\begin{equation}
Cov(I(z^S),I(z^{S'})) =\sigma^2 \frac{\lambda}{2}\frac{1}{\omega^2 n^2} \sum_{j=1}^n e^{-i\omega \theta^{(j)}
 \cdot (z^S-z^{S'})} (\theta^{(j)})^\top R_\omega(z^S,z^{S'})
 \theta^{(j)}.
\end{equation}
The signal-to-noise ratio is given by
\begin{equation}
(SNR)_I = \frac{\mathbb{E}[I(z_r)] }{(Var(I(z_r))^{\frac{1}{2}}}.
\end{equation}
If we compute
\begin{equation}
Var(I(z_r))\sim \sigma^2 \frac{\pi}{8 \omega^2 n},
\end{equation}
then $(SNR)_I$ can be expressed as
\begin{equation}
(SNR)_I=\frac{ \sqrt{ \pi n}  \delta^2 \omega^2 [\sigma_r-1]  U_I  }{ [\sigma_r +1]  \sigma }.
\end{equation}
The backpropagation functional is very stable with respect to
measurement noise. Of course, the number of measurements increases
the stability because the measurement noise is averaged out. We
will see in the following that the second-harmonic imaging is also
pretty stable with respect to measurement noise.

\subsubsection{Second-harmonic backpropagation}

We write, for each $j\in [0,n]$, $v_j$ as
\begin{equation}
v^{(j)}(x)=-\delta^2 (2\omega)^2 \left(\sum_{k,l} \chi_{k,l} \partial_{x_k}U^{(j)}(z_r) \partial_{x_l}
U^{(j)}(z_r)\right) \G{x}{z_r} + \nu^{(j)}(x),
\end{equation}
where $\nu_j$ is the measurement noise at the $j$-th measurement. Without any medium noise
the source term $(S)$ can be written as
\begin{equation}
(S)^{\theta^{(j)}}= \sum_{k,l} \chi_{k,l} \partial_{x_k}U^{(j)}(z_r) \partial_{x_l}
U^{(j)}(z_r) =- \omega^2U_I^2 e^{2i\omega \theta^{(j)} \cdot z_r} \sum_{k,l}
\chi_{k,l} \theta^{(j)}_k \theta^{(j)}_l.
\end{equation}
So we can write $J$ as
\begin{equation}
J(z^S)=\frac{1}{n}\sum_{j=1}^n \int_{ \partial \Omega}v^{(j)}(x) \overline{\G{x}{z^S}}
e^{-2i\omega \theta^{(j)} \cdot z^S} dx,
\end{equation}
or equivalently,
\begin{multline}
J(z^S)= -\delta^2 (2\omega)^2\frac{1}{n}\sum_{j=1}^n (S)^{\theta^{(j)}}\int_{\partial \Omega} \G{x}{z_r} \overline{\G{x}{z^S}}
 e^{-2i\omega \theta^{(j)} \cdot z^S} dx\\+\frac{1}{n}\sum_{j=1}^n \int_{\partial \Omega } \nu^{(j)}(x) \overline{\G{x}{z^S}}
  e^{-2i\omega \theta^{(j)} \cdot z^S}dx .
\end{multline}
We get that
\begin{equation}
\mathbb{E}[J(z^S)]= -\delta^2 (2\omega)^2\frac{1}{n}\sum_{j=1}^n (S)^{\theta^{(j)}} e^{-2i\omega \theta^{(j)}\cdot z^S} Q_{2\omega}(z_r,z^S),
\end{equation}
so that, using (\ref{EQUIVQ}):
\begin{equation}
\mathbb{E}[J(z_r)]\sim \delta^2  U_I^2 \frac{\omega^3}{2n} \sum_{k,l,j} \chi_{k,l}  \theta^{(j)}_k \theta^{(j)}_l.
\end{equation}
We can compute the covariance
\begin{multline}
Cov(J(z^S),J(z^{S'}))=\mathbb{E }\bigg[ \frac{1}{n^2}\left(\sum_{j=1}^ne^{-2i\omega \theta^{(j)} \cdot
z^S}  \int_{\partial \Omega}
\nu^{(j)}(x)  \overline{\G{x}{z^S}}  dx \right) \\ \left( \sum_{l=1}^n e^{2i\omega \theta^{(l)} \cdot
z^{S'}}  \int_{\partial \Omega}\nu^{(l)}(x)\G{x'}{z^{S'}}dx'\right) \bigg],
\end{multline}
which yields
\begin{equation}
Cov(J(z^S),J(z^{S'}))=
\sigma^2 \frac{\lambda}{2} Q_{2\omega}(z^{S'},z^S)\frac{1}{n^2}\sum_{j=1}^ne^{-2i\omega
\theta^{(j)} \cdot (z^S-z^{S'})} .
\end{equation}
Now we have
\begin{equation}
Var(J(z_r))^{1/2}\sim  \frac{\sigma}{2\omega} \sqrt{\frac{\pi}{2n}}.
\end{equation}
The signal-to-noise ratio,
\begin{equation}
(SNR)_J=\frac{\mathbb{E}[J(z_r)] }{(Var(J(z_r))^{\frac{1}{2}}},
\end{equation}
is given by
\begin{equation}
(SNR)_J= \frac{ 2\delta^2 \omega^2 U_I\left(\sum_j \sum_{k,l} \chi_{k,l}
\theta^{(j)}_k \theta^{(j)}_l \right) }{\pi\sigma\sqrt{n}}.
\end{equation}

Even though it appears that the $(SNR)$ is proportional to
$\frac{1}{\sqrt{n}}$, the term $\sum_j \theta^{(j)}_k
\theta^{(j)}_l $ is actually much bigger. In fact, if we pick
$\theta^{(j)}=\frac{2j\pi}{n}$ we get that \begin{equation}
\sum_{k,l} \chi_{k,l} \sum_j \theta^{(j)}_k \theta^{(j)}_l=
\sum_{j=1}^n \left( \chi_{1,1}\cos^2 \frac{2j\pi}{n} +
\chi_{2,2}\sin^2\frac{2j\pi}{n} +2\chi_{1,2}\sin \frac{2j\pi}{n}
\cos \frac{2j\pi}{n} \right),
\end{equation}
and hence,
\begin{equation}
\sum_{k,l} \chi_{k,l} \sum_j \theta^{(j)}_k \theta^{(j)}_l
\sim\frac{n}{2} \max[\chi_{1,1}, \chi_{2,2}] .
\end{equation}
Therefore, we can conclude that
\begin{equation}
(SNR)_J= \frac{ \delta^2 \omega^2U_I^2\sqrt{n}\max[\chi_{1,1}, \chi_{2,2}] }{\pi \sigma_\nu }.
\end{equation}
The signal-to-noise ratio is very similar to the one seen in the
classic backpropagation case. So the sensitivity with respect to
relative measurement noise should be similar. It is noteworthy
that in reality, due to very small size of the (SHG) signal
($\chi$ has a typical size of $10^{-12} \ m/V$), the measurement
noise levels will be higher for the second-harmonic signal.

\section{Numerical results} \label{sec7}
\subsection{The direct problem}
\begin{figure}[!h]
   \begin{minipage}{.45\linewidth}
   \includegraphics[width=\linewidth]{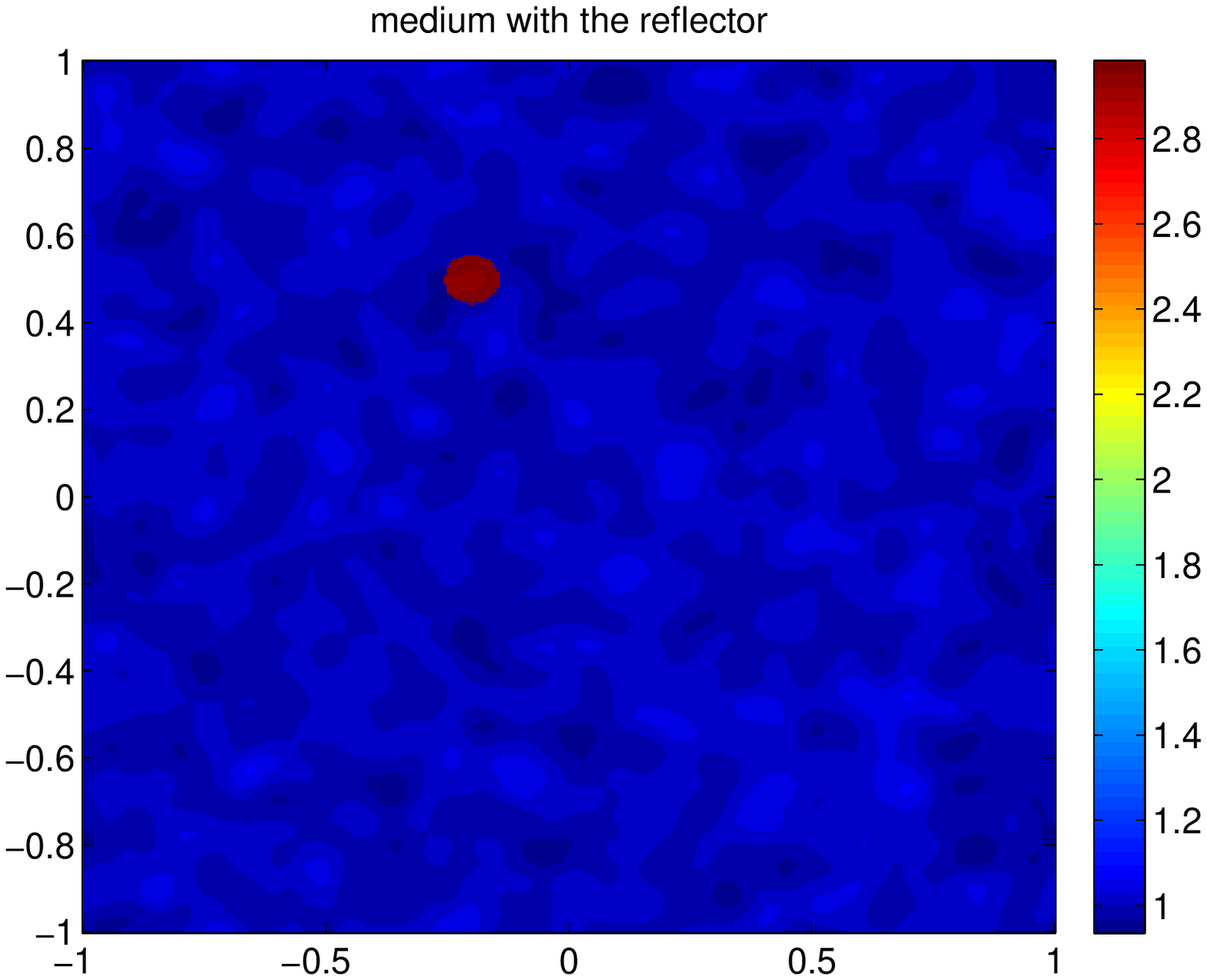}
   \caption{\label{graphM} Medium with the reflector.}
 \end{minipage}
   \begin{minipage}{.45\linewidth}
   \includegraphics[width=\linewidth]{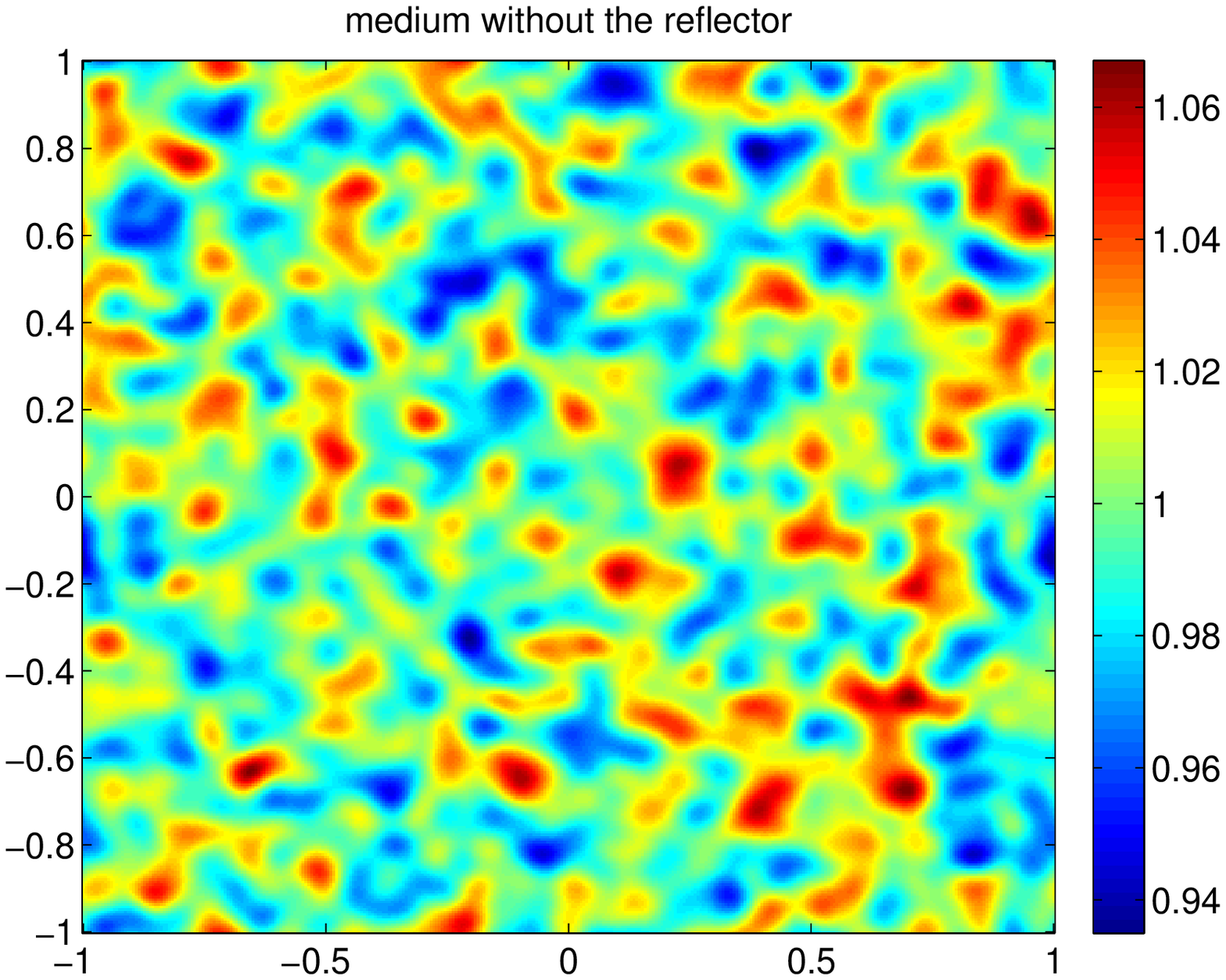}
   \caption{\label{graphMssref} Medium without the reflector (permittivity variations zoomed out).}
 \end{minipage}
\end{figure}

We consider the medium to be the square $[-1,1]^2$. The medium has
an average propagation speed of $1$, with random fluctuations with
Gaussian statistics (see Figure \ref{graphMssref}). To simulate
$\mu$ we use the algorithm described in \cite{TDerivativ} which
generates random Gaussian fields with Gaussian covariance function
and take a standard deviation equal to $0.02$ and a correlation
length equal to $0.25$. We consider a small reflector in the
medium $\Omega_r=z_r+\delta B(0,1)$ with $z_r=(-0.2,0.5)$ and
$\delta=0.004/\pi$, represented on Figure \ref{graphM}. The
contrast of the reflector is $\sigma_r=2$. We fix the frequency to
be $\omega=8$. We get the boundary data $u_{s}$ when the medium is
illuminated by the plane wave $U_{I}(x)=e^{i \omega \theta \cdot
x}$. The correlation length of the medium noise was picked so that
it has a similar size as the wavelength of the illuminating plane
wave. We get the boundary data by using an integral representation
for the field $u_{s,\theta}$. We also compute the boundary data
for the second-harmonic field $v$. We compute the imaging
functions $I$ and $J$  respectively defined in (\ref{DefI}) and
(\ref{DefJ}), averaged over two different lightning  settings.
(see Figures \ref{graphI} and \ref{graphJ} for instance).

\begin{figure}[!h]
 \begin{minipage}{.45\linewidth}
  \centering\includegraphics[width=\linewidth]{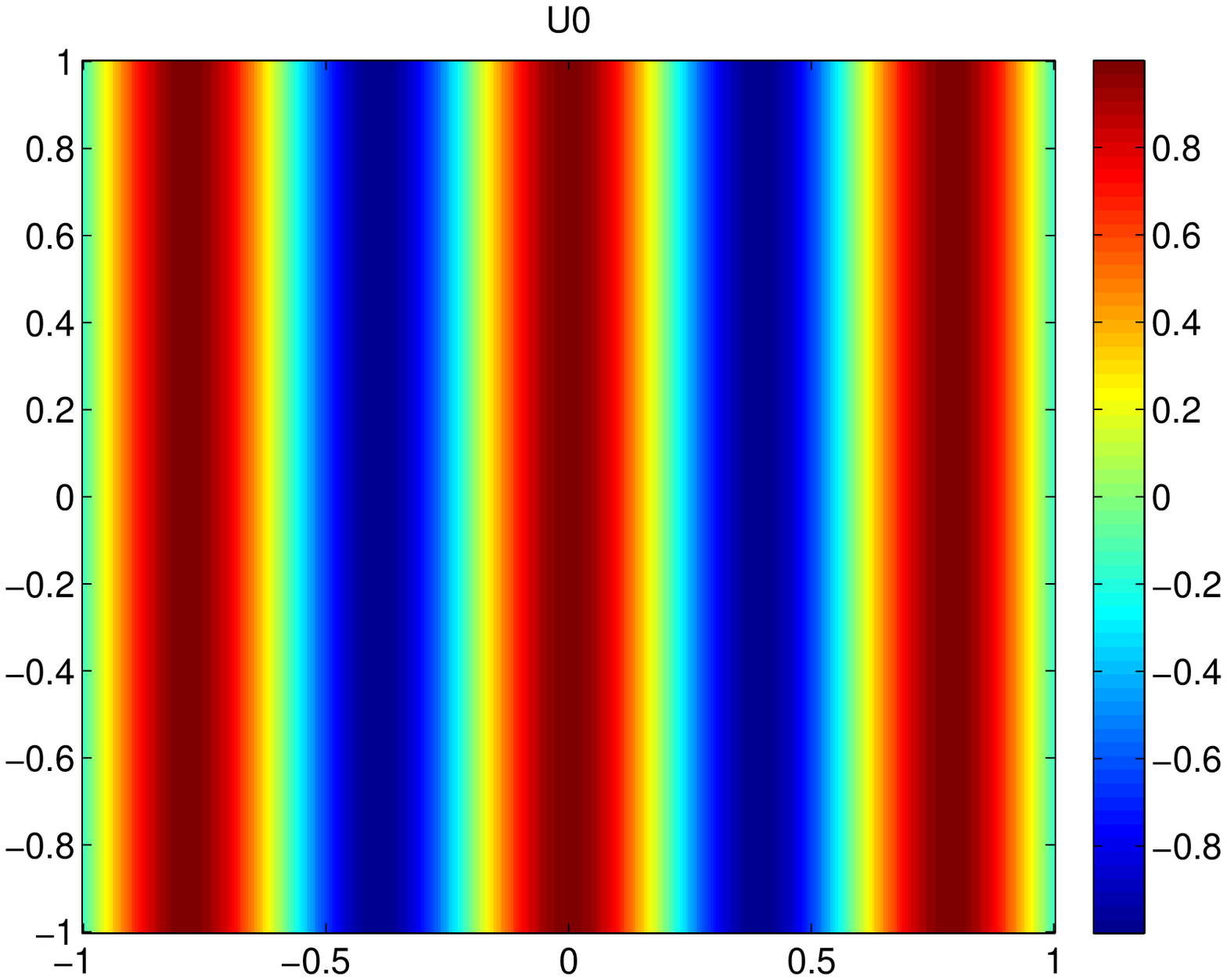}
  \caption{\label{UI} Incoming field $U_I$.}
 \end{minipage}
\begin{minipage}{.45\linewidth}
  \centering\includegraphics[width=\linewidth]{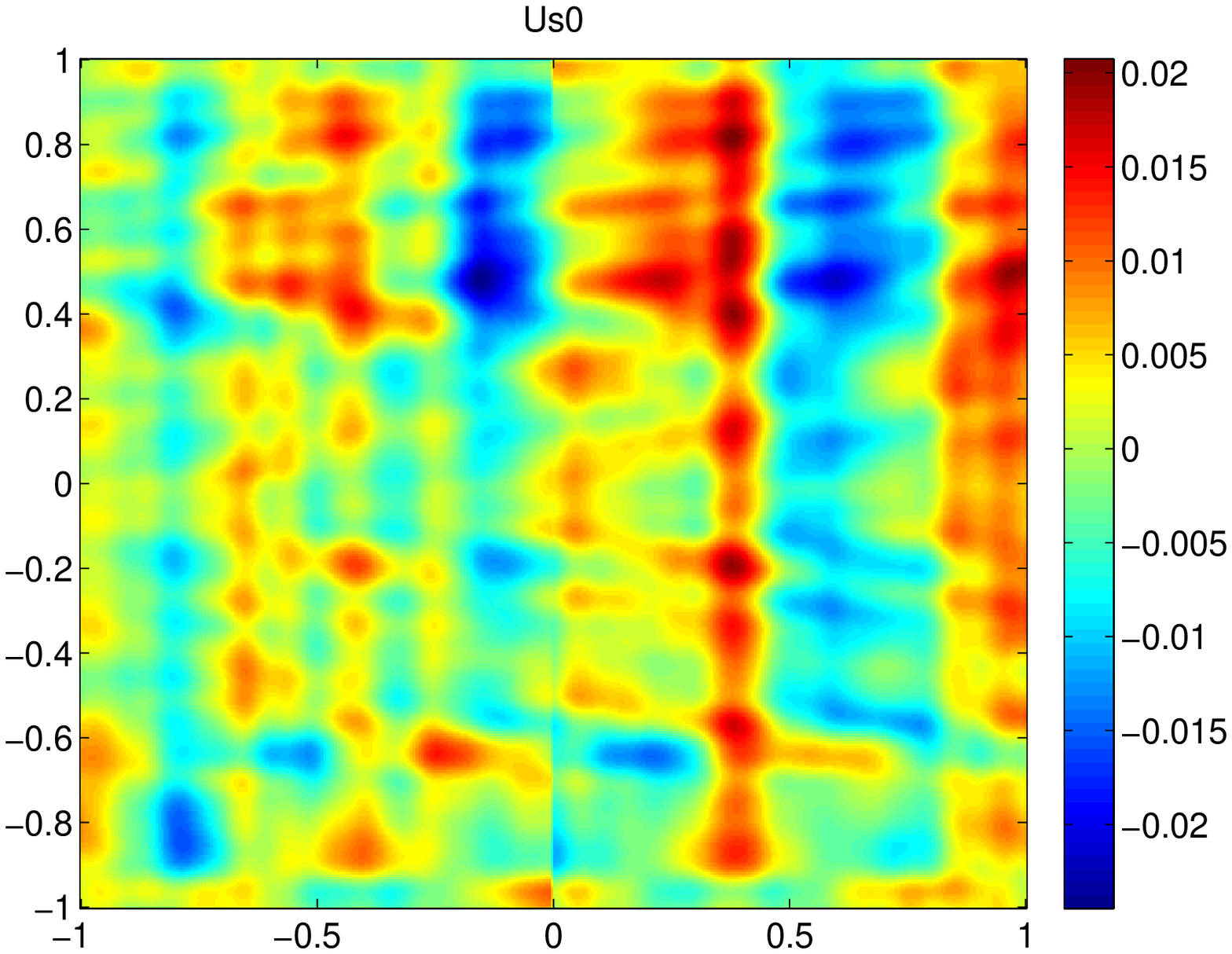}
  \caption{\label{Us0} Background field in the absence of a reflector $u_s^{(\mu)}$.}
 \end{minipage}

\end{figure}

\begin{figure}[!h]
 \begin{minipage}{.45\linewidth}
  \centering\includegraphics[width=\linewidth]{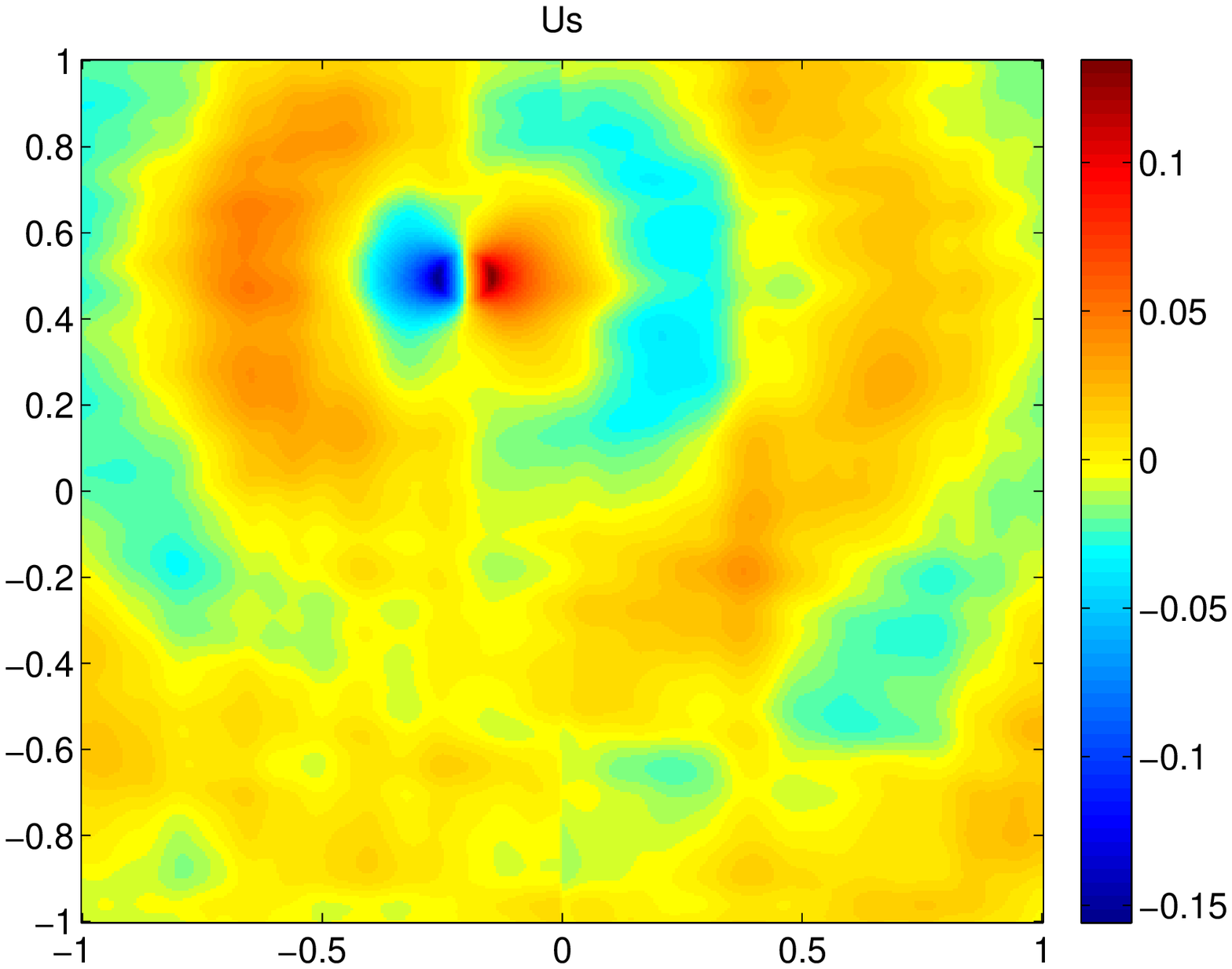}
  \caption{\label{graphus} Total scattered field $u_s$.}
 \end{minipage}
\begin{minipage}{.45\linewidth}
  \centering\includegraphics[width=\linewidth]{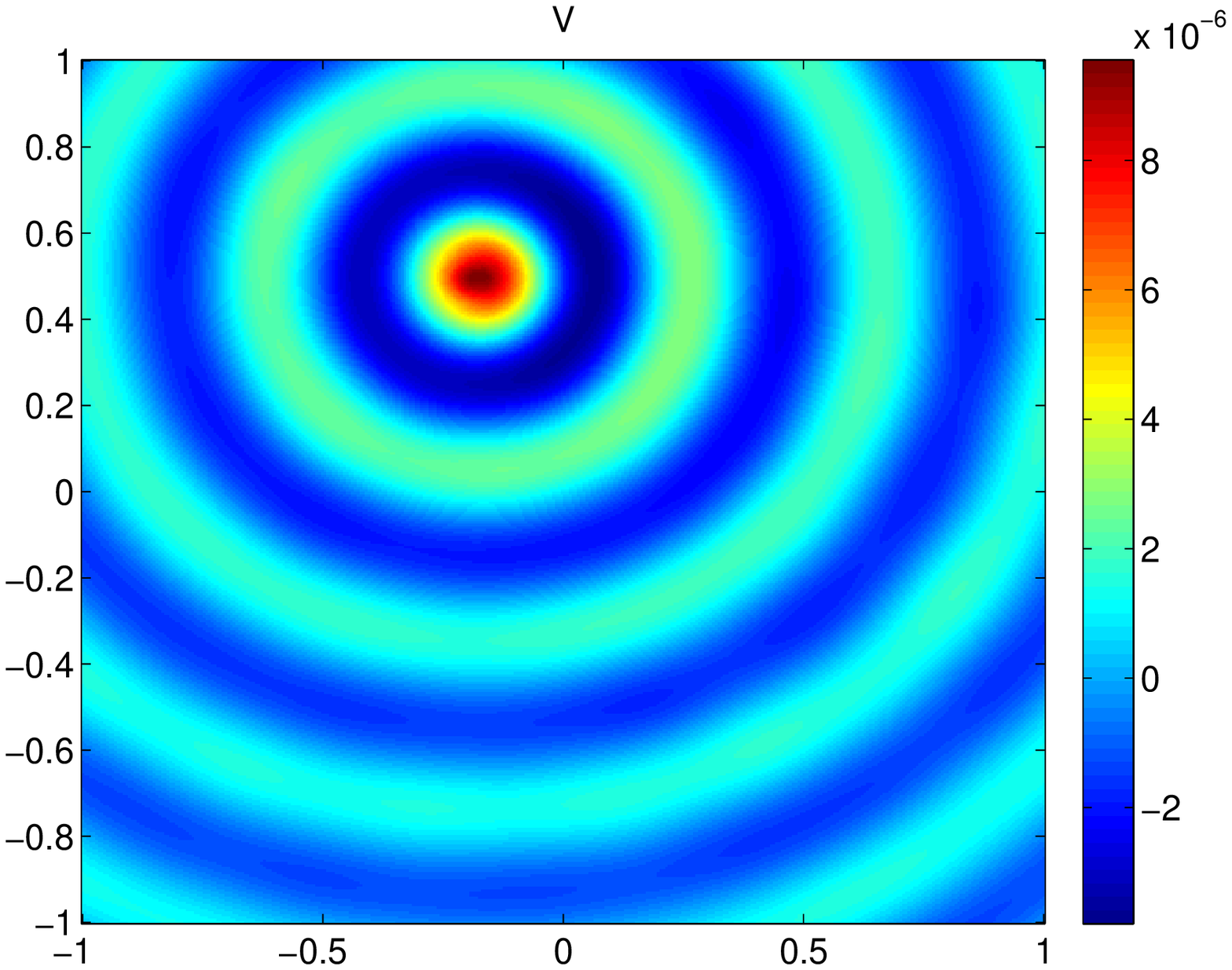}
  \caption{\label{graphItilde30} Second-harmonic field $v$.}
 \end{minipage}

\end{figure}

\subsection{The imaging functionals and the effects of the number of plane wave illuminations}
\begin{figure}[!h]
 \begin{minipage}{.45\linewidth}
  \centering\includegraphics[width=\linewidth]{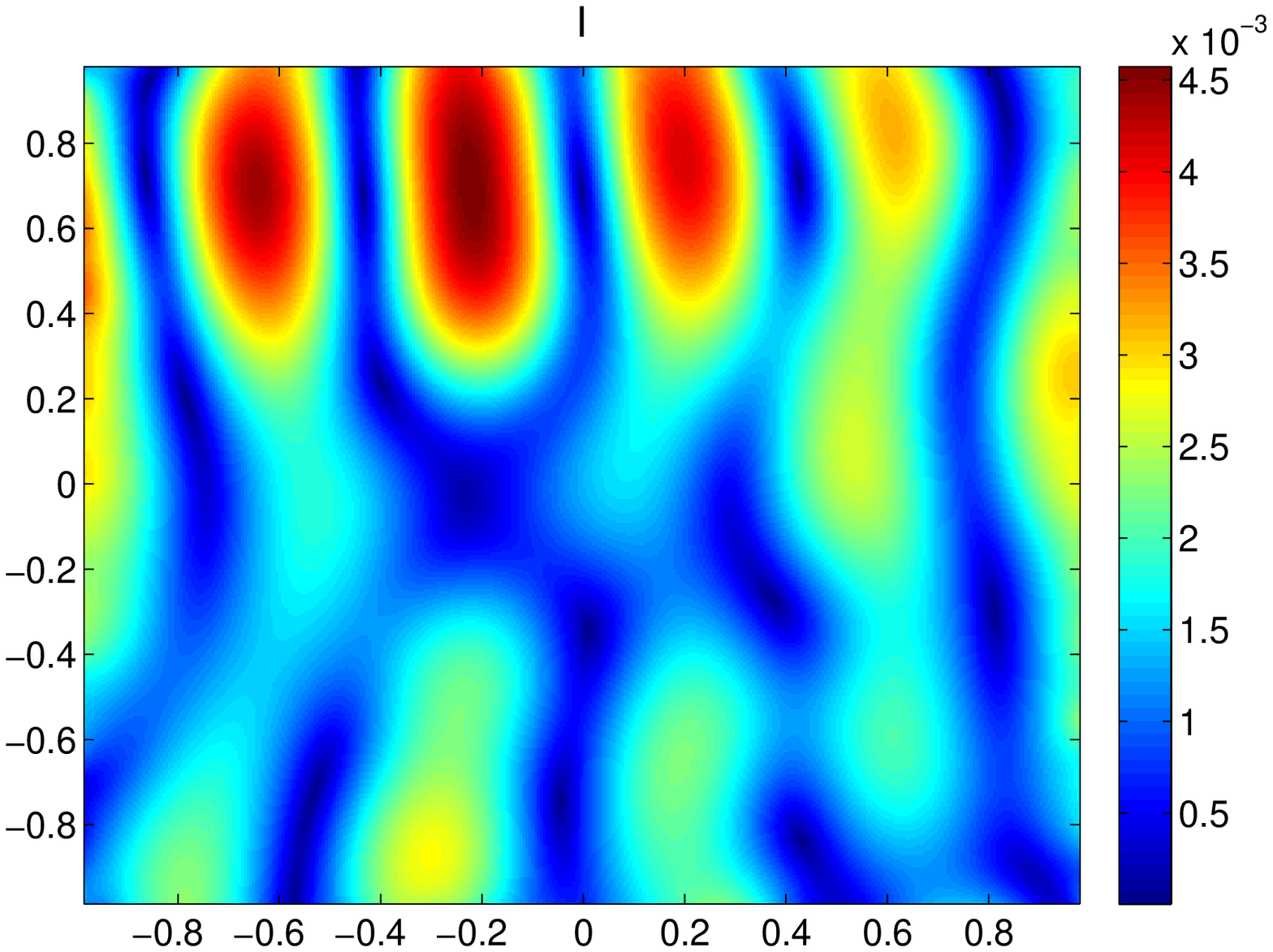}
  \caption{\label{graphI} $I$ with $1$ illumination.}
 \end{minipage}
\begin{minipage}{.45\linewidth}
  \centering\includegraphics[width=\linewidth]{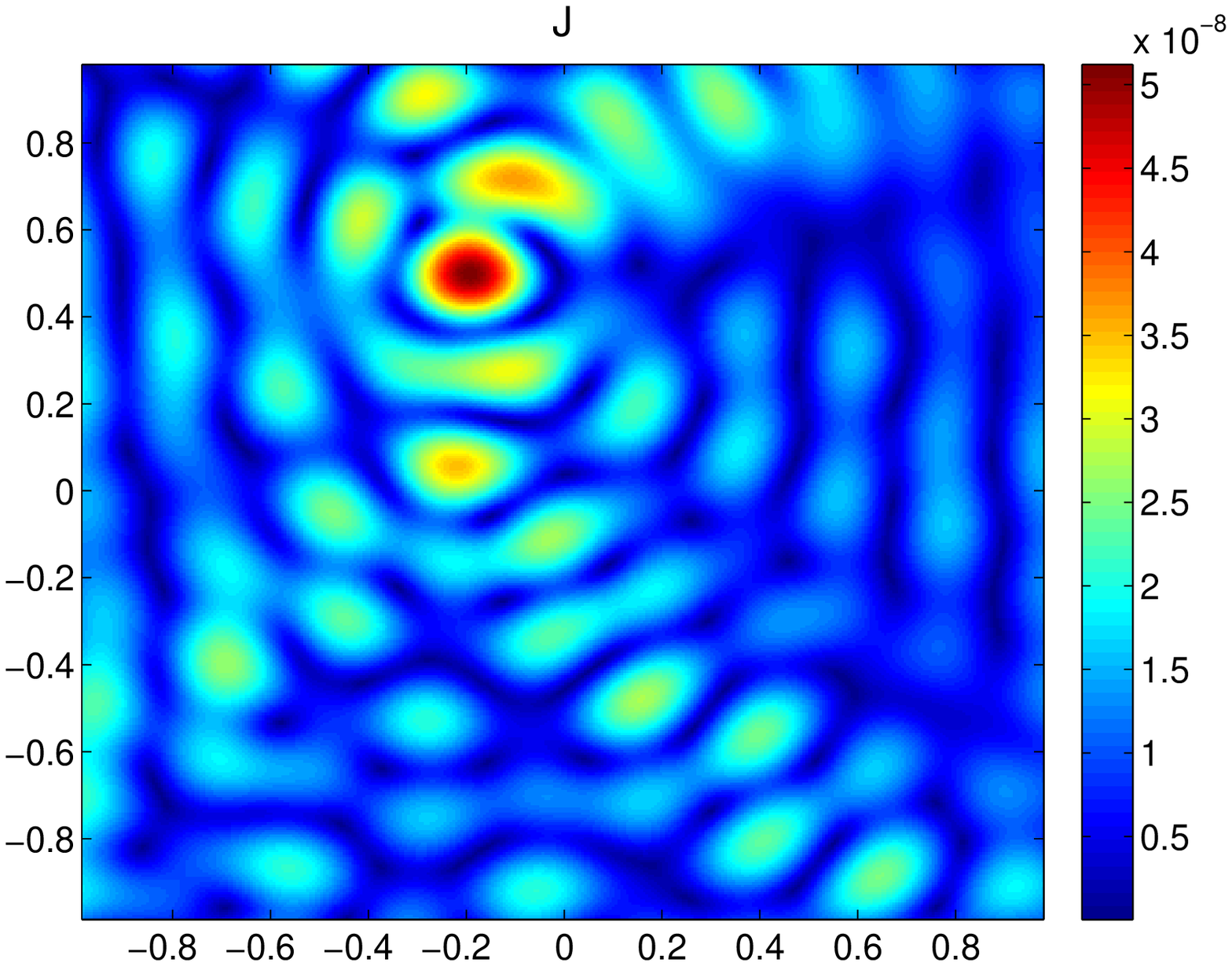}
  \caption{\label{graphJ} $J$ with $1$ illumination.}
 \end{minipage}

\end{figure}

\begin{figure}[!h]
 \begin{minipage}{.45\linewidth}
  \centering\includegraphics[width=\linewidth]{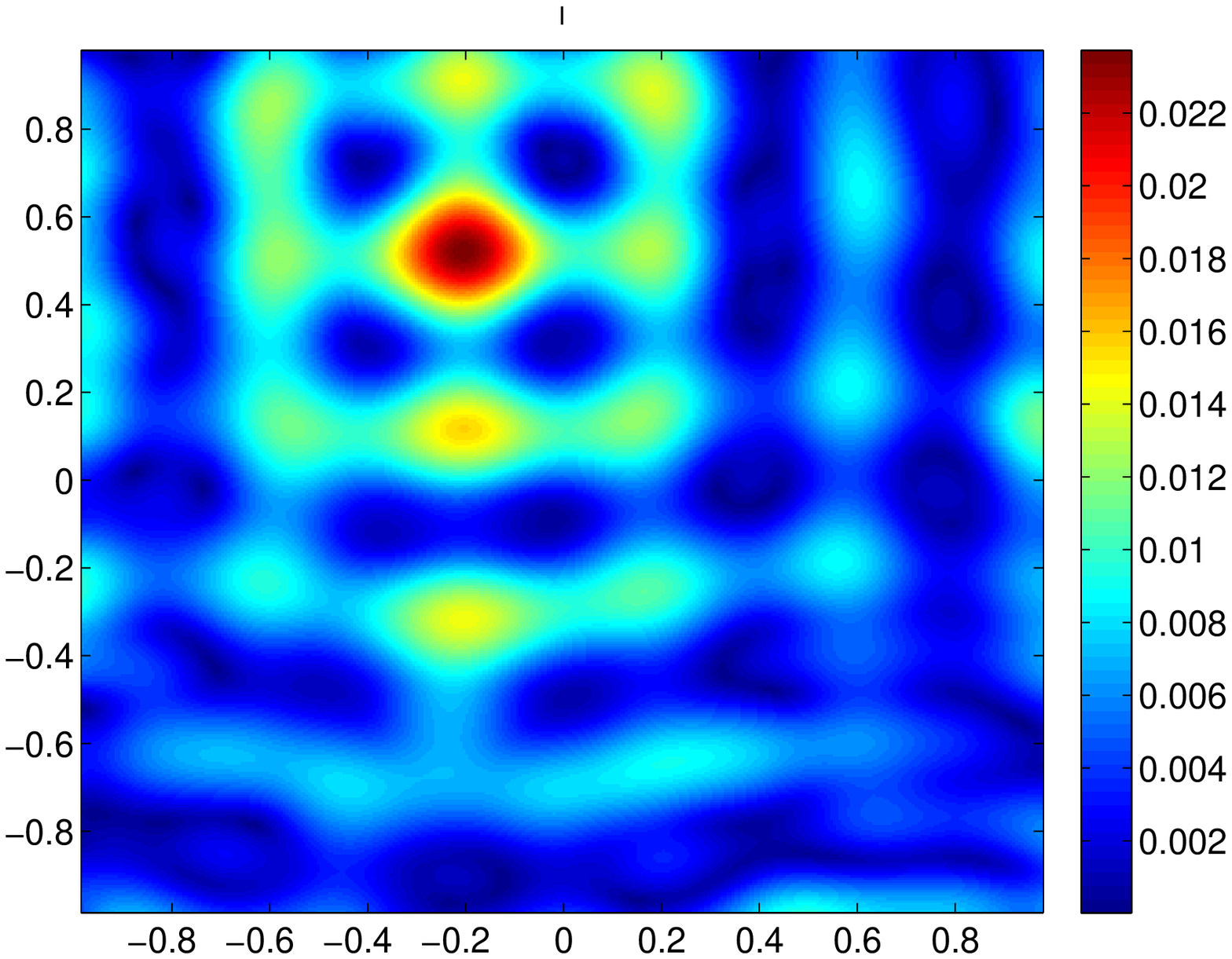}
  \caption{\label{graphI4} $I$ with $4$ illuminations.}
 \end{minipage}
\begin{minipage}{.45\linewidth}
  \centering\includegraphics[width=\linewidth]{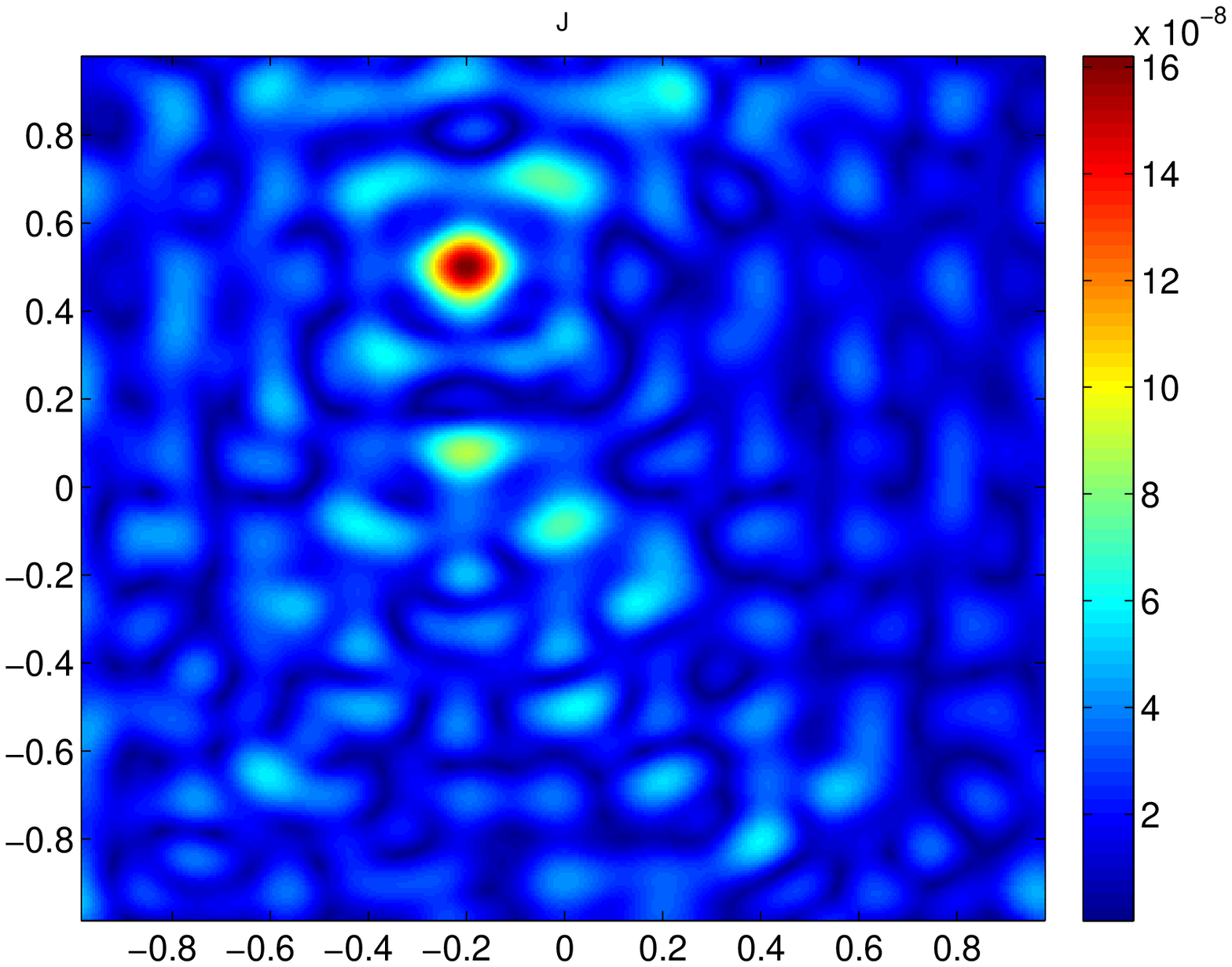}
  \caption{\label{graphJ4} $J$ with $4$ illuminations.}
 \end{minipage}

\end{figure}

\begin{figure}[!h]
 \begin{minipage}{.45\linewidth}
  \centering\includegraphics[width=\linewidth]{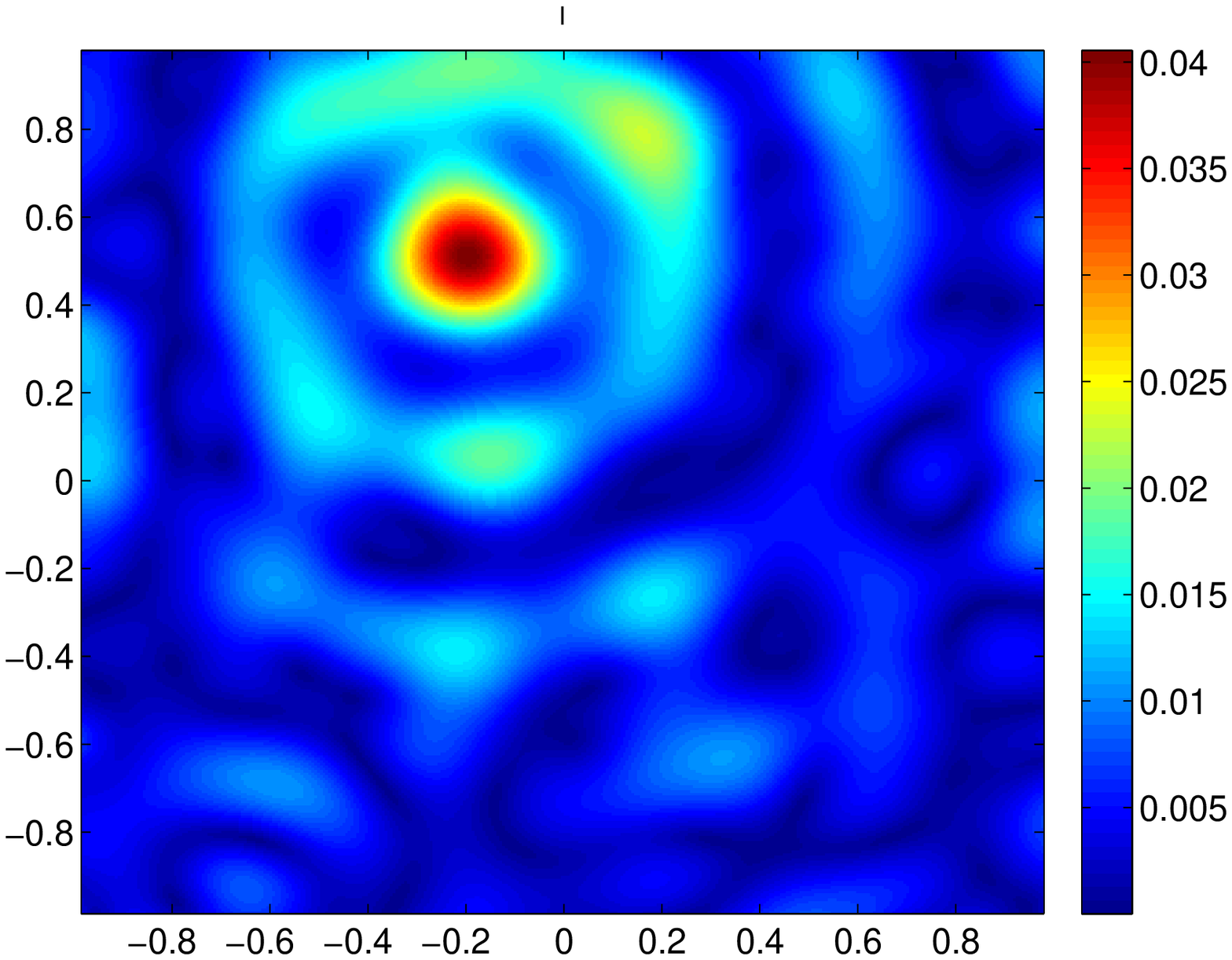}
  \caption{\label{graphI8} $I$ with $8$ illuminations.}
 \end{minipage}
\begin{minipage}{.45\linewidth}
  \centering\includegraphics[width=\linewidth]{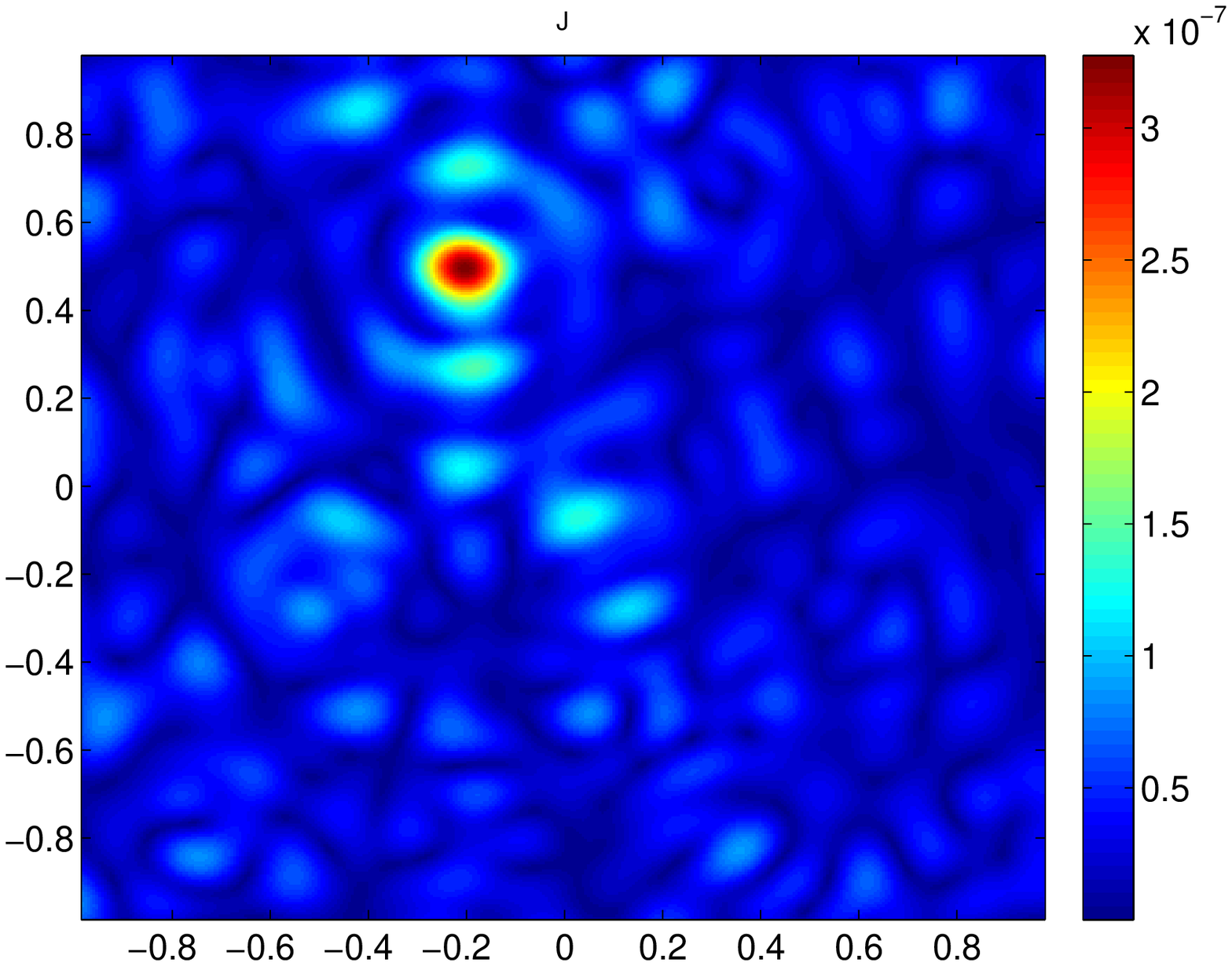}
  \caption{\label{graphJ8} $J$ with $8$ illuminations.}
 \end{minipage}

\end{figure}

\begin{figure}[!h]
 \begin{minipage}{.45\linewidth}
  \centering\includegraphics[width=\linewidth]{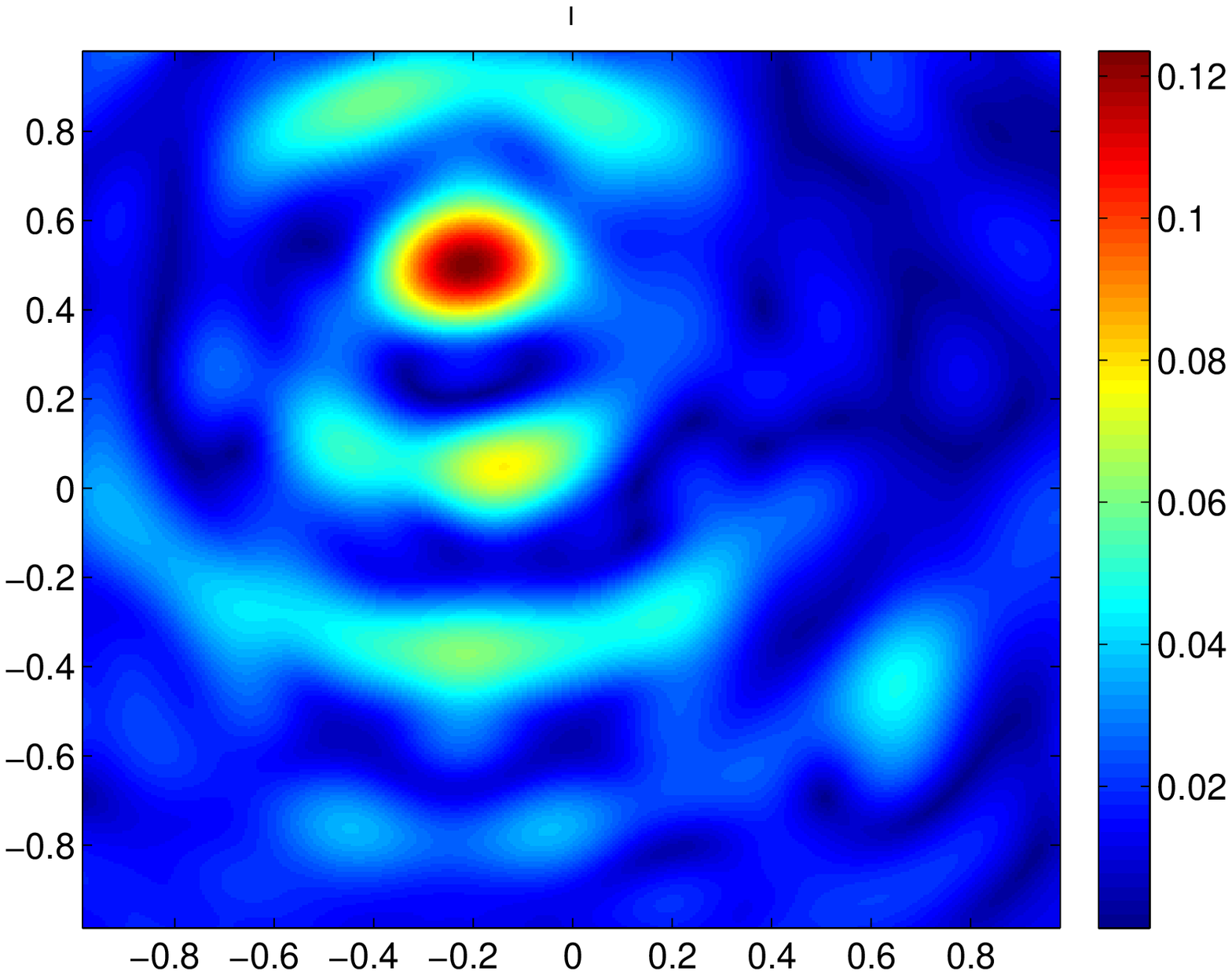}
  \caption{\label{graphI32} $I$ with $32$ illuminations.}
 \end{minipage}
\begin{minipage}{.45\linewidth}
  \centering\includegraphics[width=\linewidth]{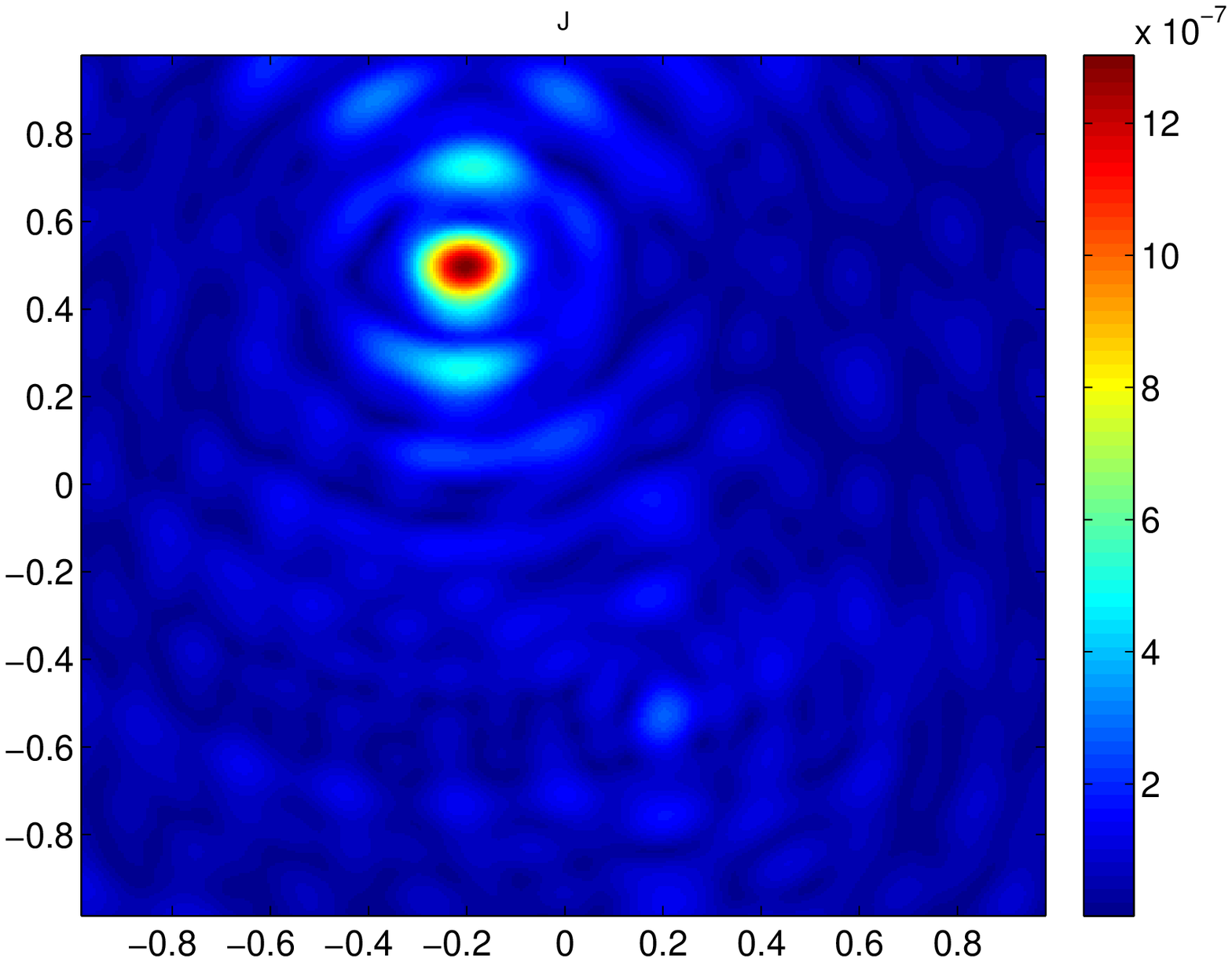}
  \caption{\label{graphJ32} $J$ with $32$ illuminations.}
 \end{minipage}

\end{figure}

We compute the imaging functionals $I$ and $J$  respectively
defined in (\ref{DefI}) and (\ref{DefJ}), averaged over four
different illuminations settings. We fix the noise level
($\sigma_\mu =0,02$), the volume of the particle ($v_r=10^{-2}$)
and the contrast $\sigma_r = 2$. In Figures~\ref{graphI} and
\ref{graphJ} the image is obtained after backpropagating the
boundary data from one illumination ($\theta =0$). On the
following graphs, we average over several illumination
angles:\begin{itemize} \item $4$ uniformly distributed angles for
Figures~\ref{graphI4} and ~\ref{graphJ4}. \item $8$ uniformly
distributed angles for Figures~\ref{graphI8} and ~\ref{graphJ8}.
\item $32$ uniformly distributed angles for
Figures~\ref{graphI32} and ~\ref{graphJ32}.
\end{itemize}
As predicted, the shape of the spot on the fundamental frequency
imaging is very dependant on the illumination angles, whereas with
second-harmonic imaging we get an acceptable image with only one
illumination. In applications, averaging over different
illumination is useful because it increases the stability with
respect to measurement noise. It is noteworthy that, as expected,
the resolution of the second-harmonic image is twice higher than
the regular imaging one.

\subsection{Statistical analysis}
\subsubsection{Stability with respect to medium noise}

\begin{figure}[!h]
\centering\includegraphics[width=\linewidth]{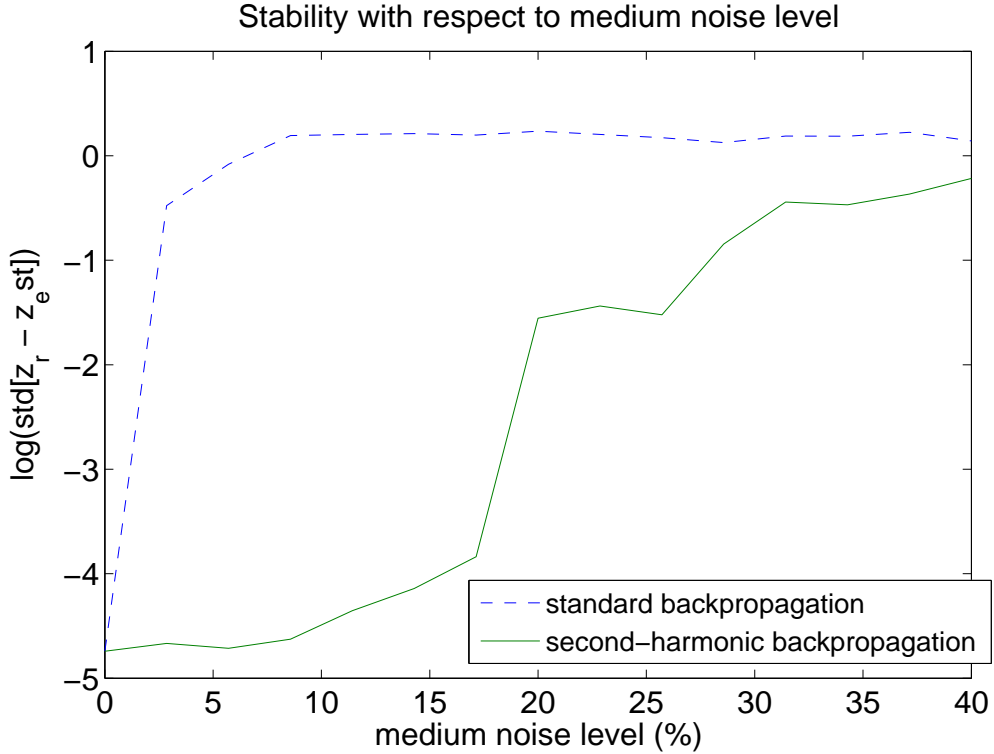}
\caption{\label{graphbruitIJ} Standard deviation of the
localization error with respect to the medium noise level for
standard backpropagation (top) and second-harmonic image
(bottom).}
\end{figure}

Here we show numerically that the second-harmonic imaging is more
stable with respect to medium noise. In Figure~\ref{graphbruitIJ},
we plot the standard deviation of the error $\vert
z_{est}-z_r\vert $ where $z_{est}$ is the estimated location of
the reflector. For each level of medium noise we compute the error
over $120$ realizations of the medium, using the same parameters,
as above. The functional imaging $J$ is clearly more robust than
earlier.

\subsubsection{Effect of the volume of the particle}

\begin{figure}[!h]
\centering\includegraphics[width=\linewidth]{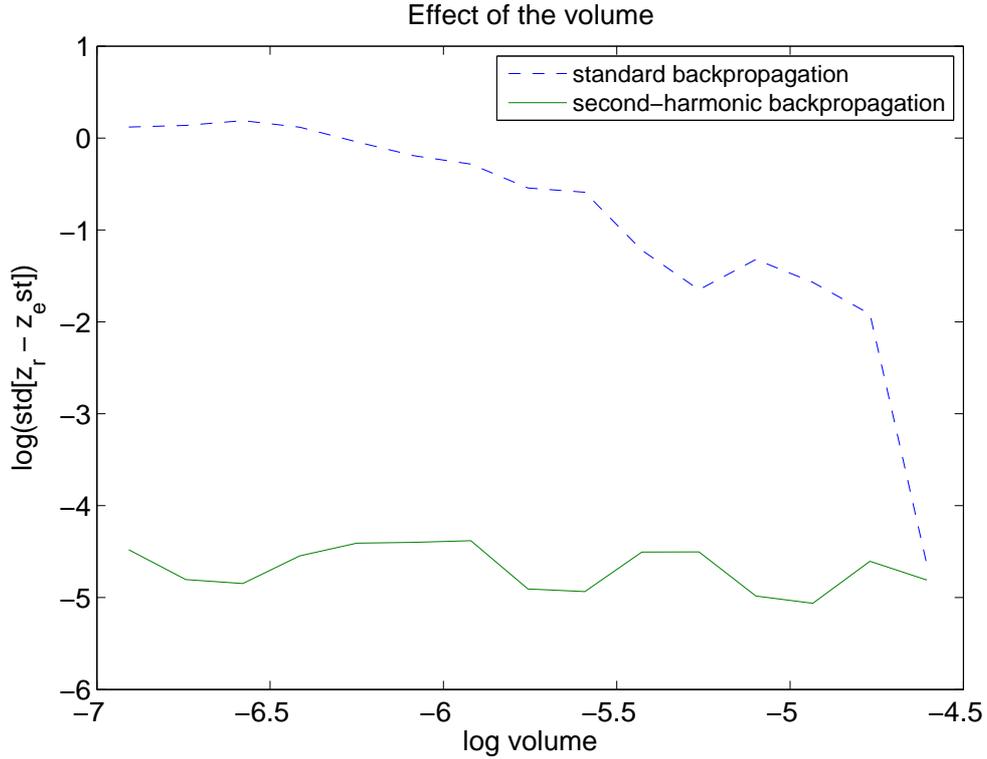}
\caption{\label{graphvolumeIJ} Standard deviation of the
localization error with respect to the reflector's volume (log
scale) for standard backpropagation (top) and second-harmonic
image (bottom).}
\end{figure}

We show numerically that the quality of the second-harmonic image
does not depend on the volume of the particle. We fix the medium
noise level ($\sigma_\mu =0.02$) and plot the standard deviation
of the error with respect to the volume of the particle
(Figure~\ref{graphvolumeIJ}). We can see that if the particle is
too small, the fundamental backpropagation algorithm cannot
differentiate the reflector from the medium and the main peak gets
buried in the speckle field. The volume of the particle does not
have much influence on the second-harmonic image quality.

\subsubsection{Stability with respect to measurement noise}

\begin{figure}[!h]
\centering\includegraphics[width=\linewidth]{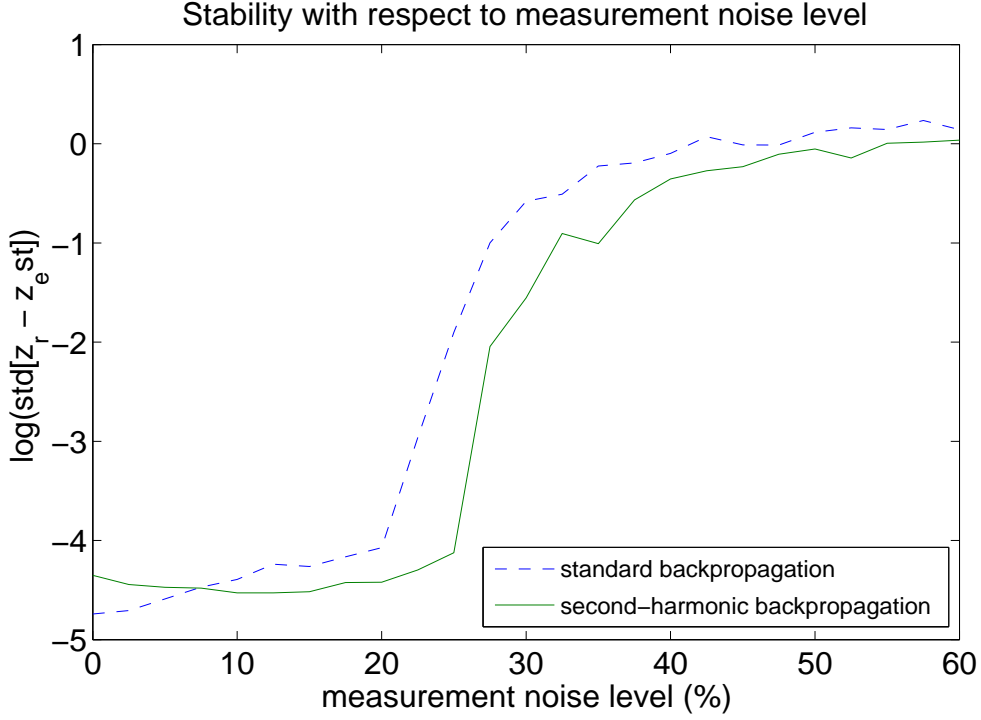}
\caption{\label{graphmeasnoise} Standard deviation of the
localization error with respect to measurement noise level for
standard backpropagation (top) and second-harmonic image
(bottom).}
\end{figure}

We compute the imaging functionals with a set of data
 obtained without any medium noise and perturbed with
a Gaussian white noise for each of $8$ different illuminations.
For each noise level, we average the results over $100$ images.
Figure~\ref{graphmeasnoise} shows that both functionals have
similar behaviors.

\begin{figure}[!h]
\begin{center}
\includegraphics[width=\linewidth]{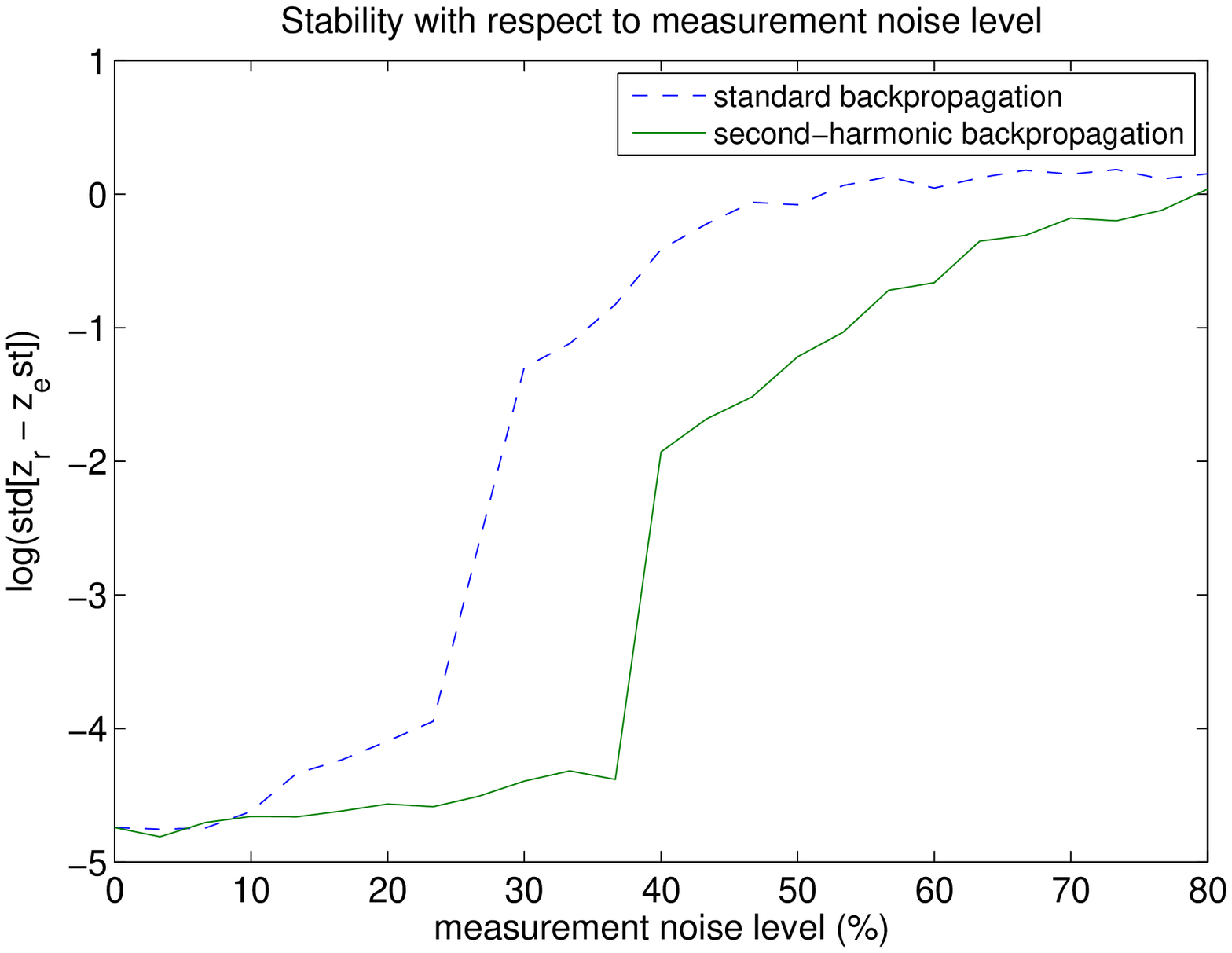}
\caption{\label{graphmeasnoise10} Standard deviation of the localization error with respect to measurement
noise level, when averaged over $16$ illuminations of angles uniformly distributed between $0$ and
$2\pi$ for standard backpropagation (top) and second-harmonic image (bottom).}
\end{center}
\end{figure}

As mentioned before, in applications, the weakness of the SHG
signal will induce a much higher relative measurement noise than
in the fundamental data. Since the model we use for measurement
noise has a zero expectation, averaging measurements over
different illuminations can improve the stability significantly as
shown in Figure~\ref{graphmeasnoise10}, where the images have been
obtained with $16$ illuminations instead of $8$.

\section{Concluding remarks}
We have studied how second-harmonic imaging can be used to locate
a small reflector in a noisy medium, gave asymptotic formulas for
the second-harmonic field, and investigated statistically the
behavior of the classic and second-harmonic backpropagation
functionals. We have proved that the backpropagation algorithm is
more stable with respect to medium noise. Our results can also be
extended to the case of multiple scatterers as long as they are
well-separated.

\appendix
\section{Proof of (\ref{estimateborn})}  \label{appenda}
Let $R$ be large enough so that $\Omega_\mu \Subset B_R$, where
$B_R$ is the ball of radius $R$ and center $0$. Let $S_R =\partial
B_R$ be the sphere of radius $R$, and introduce the
Dirichlet-to-Neumann operator $\mathcal{T}$ on $S_R$:
 \begin{equation}
\begin{aligned}\label{DefDTN}
\mathcal{T} \ : \ H^{1/2}(S_R)& \longrightarrow  H^{-1/2}(S_R) \\
u  &\longmapsto \mathcal{T}[ u].
\end{aligned}
\end{equation}
According to \cite{nedelec92}, $\mathcal{T}$ is continuous and
satisfies
\begin{equation}\label{ESTT1}
-\text{Re }\left< \mathcal{T} [u], u \right> \geq \frac{1}{2R}
\Vert u \Vert_{L^2(S_R)}^2,\quad \forall u \in H^{1/2}(S_R),
\end{equation} and
\begin{equation}\label{ESTT2}
\text{Im } \left< \mathcal{T} [u], u \right> >0 \ \text{if } u
\neq 0.
\end{equation}
Here, $\left< \,, \,\right> >$ denotes the duality pair between
$H^{1/2}(S_R)$ and $H^{-1/2}(S_R)$.  Now introduce the continuous
bilinear form $a$:
 \begin{equation}\label{DEFT}
\begin{aligned}
H^{1}(B_R)\times H^1(B_R) & \longrightarrow  \mathbb{C} \\
(u,v)  &\longmapsto a(u,v)= \int_{B_R} (1+\mu) \nabla u \cdot \overline{\nabla v} -
\omega^2 \int_{B_R} u\overline{v} -  \left<\mathcal{T}[u],v \right>,
\end{aligned}
\end{equation}
as well as the continuous bilinear form $b$:
\begin{equation}
\begin{aligned}
H^{1}(B_R)& \longrightarrow  \mathbb{C} \\
v  &\longmapsto b(v) = \int_{B_R} \mu \nabla U_0 \cdot \overline{\nabla v}.
\end{aligned}
\end{equation}
Problem (\ref{sgheq1}) has the following variational formulation:
Find $u\in H^1(B_R)$ such that
\begin{equation} \label{variationnal}
a(u,v)=b(v) \ \  \forall v \in H^1(B_R).
\end{equation}
With (\ref{ESTT1}) one can show that
\begin{equation} \label{coerca}
\text{Re } a(u,u) \geq C_1 \Vert \nabla u \Vert_{L^2(B_R)}^2 - C_2
\Vert u \Vert_{L^2(B_R)}^2,
\end{equation} so that $a$ is weakly coercive with respect to the pair
 $\left(H^1(B_R), L^2(B_R)\right)$. Since the imbedding  of $H^1(B_R)$ into
 $L^2(B_R)$ is compact we can apply Fredholm's alternative to problem
 (\ref{variationnal}). Hence, we deduce existence of a solution from uniqueness of a solution which
 easily follows by using identity (\ref{ESTT2}).

Now we want to prove that if $u$ is the solution of (\ref{variationnal}) then
\begin{equation}
\Vert u \Vert_{H^1(B_R)} \leq  \Vert \mu \Vert_\infty.
\end{equation}
We proceed by contradiction. Assume that $\forall n \in
\mathbb{N}$, there exists $\mu_n \in L^\infty (B_R)$  compactly
supported and $ u_n \in H^1(B_R)$ solution of (\ref{variationnal})
such that
\begin{equation}\label{defun}
\Vert u_n \Vert_{H^1(B_R)} \geq nC \Vert \mu_n \Vert_\infty.
\end{equation} Consider the sequence: \begin{equation} v_n=\frac{u_n}{\Vert u_n \Vert_{H^1(B_R)}}. \end{equation} $(v_n)_{n\in \mathbb{N}}$ is bounded in $H^1(B_R)$ so there exists a subsequence still denoted by $v_n$ and $v^* \in H^1(B_R)$
such that $v_n \rightharpoonup v^*$ in $H^1(B_R)$ and $v_n
\rightarrow v^*$ in $L^2(B_R)$. Now since $u_n$ is a solution of
(\ref{variationnal}), we have
\begin{equation}
\int_{B_R} (1+\mu_n) \nabla v_n \cdot \overline{\nabla v_n}  - \omega^2 \int_{B_R} v_n \overline{v_n} - \left< \mathcal{T} v_n, v_n \right> = \int_{B_R}\mu_n \nabla U_0 \cdot \overline{\nabla v_n}.
\end{equation}
Using (\ref{defun}) we obtain that
\begin{equation}
\int_{B_R} (1+\mu_n) \vert \nabla v_n \vert^2  - \omega^2 \int_{B_R} \vert v_n \vert^2  - \left< \mathcal{T} v_n, v_n \right> \longrightarrow 0 \ (n\rightarrow \infty).
\end{equation} Since  $\int_{B_R} \mu_n \vert \nabla v_n \vert^2 \longrightarrow 0$, we get that $\widetilde{a}(v_n,v_n) \longrightarrow 0$, where
\begin{equation}
\widetilde{a}(u,v)=\int_{B_R} \nabla u \cdot \overline{\nabla v}  - \omega^2 \int_{B_R} u\overline{v}  - \left< \mathcal{T} u,v\right>.
\end{equation}
We want to prove that $v_n$ converges strongly in $H^1(B_R)$ to $v^*$ and that $v^*=0$. This will
contradict the fact that $\forall n, \ \Vert v_n \Vert_{H^1(B_R)}=1$.

Now we decompose $\widetilde{a}=\widetilde{a_c}+\widetilde{a_w}$ into a coercive part \begin{equation}
\widetilde{a_c}(u,v) = \int_{B_R} \nabla u \cdot \overline{\nabla v} -\left< \mathcal{T} u,v\right>
\end{equation} and a weakly continuous part:
\begin{equation}
\widetilde{a_w}(u,v) = - \omega^2 \int_{B_R} u\overline{v}.
\end{equation}
So $\widetilde{a}(v_n-v^*,v_n-v^*)=\widetilde{a_c}(v_n-v^*,v_n-v^*) + \widetilde{a_w}(v_n-v^*,v_n-v^*)$. We write $\widetilde{a_c}(v_n-v^*, v_n-v^*)= \widetilde{a_c}(v_n-v^*) - \widetilde{a_c}(v_n -v^*,v^*)$. Now, since $v_n \rightharpoonup v$ in $H^1(B_R)$ and $\widetilde{a_c}$ is  strongly continuous on $H^1(B_R)^2$ we have that $\widetilde{a_c}(v_n-v^*,v^*) \longrightarrow 0$, and $\widetilde{a_c}(v_n-v^*,v_n)=\widetilde{a_c}(v_n,v_n)-\widetilde{a_c}(v^*,v_n) \longrightarrow -\widetilde{a_c}(v^*,v^*)$ which is
\begin{equation}
\widetilde{a_c}(v_n-v^*,v_n-v^*) \longrightarrow -\widetilde{a_c}(v^*,v^*).
\end{equation} The coercivity of $\widetilde{a_c}$ gives \begin{equation}
\widetilde{a_c}(v^*,v^*)=0
\end{equation} By a computation similar to the one just above, we also find that
\begin{equation}
\widetilde{a}(v_n-v^*,v_n-v^*) \longrightarrow -\widetilde{a}(v^*,v^*).
\end{equation}
Since $\widetilde{a_w}(v_n-v^*,v_n-v^*)\longrightarrow 0$, we get that \begin{equation}
\widetilde{a}(v^*,v^*)=0.
\end{equation} So $v^*$ =0 and, since $\widetilde{a}$ satisfies (\ref{coerca}), we get that $\Vert \nabla v_n \Vert^2_{L^2(B_R)} \longrightarrow 0$ as $n\rightarrow \infty$. We have
\begin{equation}
v_n \longrightarrow v=0 \text{ in } H^1(B_R).
\end{equation}

\section{Proof of Proposition \ref{propb}} \label{appendb}
Denote $V=u_s- u_s^{(\mu)}-w^{(\mu)}\cdot \nabla U_0(z_r)$. $V$ is a solution on $\mathbb{R}^2$ of
\begin{equation}
\nabla \cdot (1+\mu + [\sigma_r -1]\textbf{1}_{\Omega_r}) \nabla V + \omega^2 V= -\nabla \cdot
[\sigma_r -1]\textbf{1}_{\Omega_r}  \nabla \left[ U_0- \nabla (x-z_r) \cdot \nabla U_0(z_r)\right]
\end{equation} subject to the Sommerfeld radiation condition.
Now, define $V_0$ as the solution on $\mathbb{R}^2$ of:
\begin{equation}
\nabla \cdot (1+\mu + [\sigma_r -1]\textbf{1}_{\Omega_r}) \nabla V_0 = -\nabla \cdot [\sigma_r -1]\textbf{1}_{\Omega_r}  \nabla \left[ U_0- \nabla (x-z_r) \cdot \nabla U_0(z_r)\right].
\end{equation} with the condition $V_0(x) \longrightarrow 0 \ (x\rightarrow \infty) $. %\frac{\partial V_0}{\partial \nu}(x) =0$, $\forall x\in \partial B_R$.

From \cite[Lemma A.1]{impedio}, there exist three positive
constants $C$, $\widetilde{C}$ and $\kappa$ independent of $\mu$
and $\delta$ such that
\begin{equation}\label{appendBestV1}
\Vert \nabla V_0 \Vert_{L^2(B_R)} \leq C \delta \Vert \nabla \left[ U_0- \nabla (x-z_r) \cdot \nabla U_0(z_r)\right] \Vert_{L^\infty (\Omega_r)},
\end{equation} and
\begin{equation}\label{appendBestV2}
\Vert V_0\Vert_{L^2(B_R)} \leq \widetilde{C} \delta^{1+\kappa} \Vert \nabla \left[ U_0- \nabla (x-z_r) \cdot \nabla U_0(z_r)\right] \Vert_{L^\infty (\Omega_r)}.
\end{equation}
If we write $W=V-V_0$, we have that $W$ solves:
\begin{equation}\label{eqW}
\nabla \cdot (1+\mu + [\sigma_r -1]\textbf{1}_{\Omega_r}) \nabla W +\omega^2 W  = - \omega^2 V_0,
\end{equation} with the boundary condition $\frac{\partial W}{\partial \nu}-\mathcal{T}_\omega (W)=
 \mathcal{T}_\omega(V) - \mathcal{T}_0 (V_0)$ on $\partial B_R$, where $\mathcal{T}_\omega$ is the
 Dirichlet-to-Neumann map on $S_R$ defined in (\ref{DefDTN}) associated with the frequency $\omega$.
 The condition can be re-written : $\frac{\partial W}{\partial \nu}-\mathcal{T}_\omega (W)=
 \left(\mathcal{T}_\omega - \mathcal{T}_0\right) (V_0)$. So, based on the well posedness of (\ref{eqW}),
 there exist a constant $C'$ independent of $\mu$ and $\delta$ such that
\begin{equation}
\Vert W \Vert_{H^1(B_R)} \leq C' \left( \Vert V_0 \Vert_{L^2(B_R)} + \Vert \left[\mathcal{T}_\omega - \mathcal{T}_0         \right](V_0)\Vert_{L^2(\partial B)} \right) .
\end{equation} Now, we can write that, for some constant still denoted $C$ independent of $\mu$ and $\delta$:
\begin{equation}
\Vert V \Vert_{H^1(B_R)} \leq C \left(\Vert V_0 \Vert_{H^1(B_R)}+\Vert V_0 \Vert_{L^2(B_R)}\right).
\end{equation} Since $\delta <1$, using (\ref{appendBestV1}) and (\ref{appendBestV2}) we get
\begin{equation}
\Vert V \Vert_{H^1(B_R)} \leq C \delta^2.
\end{equation}

\section{Proof of Proposition \ref{propappendc}} \label{appendc}

Denote $\phi$: $y\longrightarrow \widetilde{y}=\phi(y)=\frac{y-z_r}{\delta}$. If we define
$\forall \widetilde{y} \in B(0,1)$:   $\widetilde{w}^{(\mu)}(\widetilde{y})=\frac{1}{\delta}
 w^{(\mu)}(\phi^{-1}(\widetilde{y}))$, we want to prove the following:
\begin{equation}
\Vert \widetilde{w}^{(\mu)}(\widetilde{y}) -\widetilde{y} -\widetilde{w}(\widetilde{y}) \Vert_{H^1(B(0,1))}
\leq C \left( \Vert \mu \Vert_{\infty} + \delta \omega^2\right).
\end{equation}
Now, using  (\ref{eqwmu}), one can see that $\widetilde{w}^{(\mu)}$ satisfies the following equation:
\begin{equation}
\nabla \cdot \left( 1+ [\sigma_r-1] \mathbf{1}_B + \widetilde{\mu} \right) \nabla \widetilde{w}^{(\mu)} + \omega^2 \delta \widetilde{w}^{(\mu)} = \nabla \cdot \left( [\sigma_r -1] \mathbf{1}_B \nabla \widetilde{y} \right),
\end{equation}
where $\widetilde{\mu}= \mu \circ \phi^{-1}$, equipped with the
Sommerfeld radiation condition. Using equation (\ref{fv1conduc})
we get that
\begin{equation}
\nabla \cdot \left( 1+ [\sigma_r-1] \mathbf{1}_B + \widetilde{\mu} \right) \nabla\left( \widetilde{w}^{(\mu)} -\widetilde{y} - \widetilde{w}\right) = - \nabla \cdot \left(\widetilde{\mu} \nabla \widetilde{w}^{(\mu)}\right) - \omega^2 \delta \widetilde{w}^{(\mu)},
\end{equation} % with the boundary condition on $B(0,1)$ :
%\begin{equation}
%\frac{\partial}{\partial \nu} \left(\widetilde{w}^{(\mu)} -\widetilde{y} - \widetilde{w}\right) -\mathcal{T}_0 %\left(\widetilde{w}^{(\mu)} -\widetilde{y} - \widetilde{w} \right)=[\mathcal{T}_{\delta \omega^2} - \mathcal{T}_0 ] %\widetilde{w}^{(\mu)}
%\end{equation}
 Now, using Meyer's theorem \cite{meyers},  we get the following estimate:
\begin{equation}\label{estdiffwtild}
\Vert \nabla\left(\widetilde{w}^{(\mu)}(\widetilde{y}) -\widetilde{y} -\widetilde{w}(\widetilde{y})\right) \Vert_{L^2(B)} \leq C \left( \Vert \widetilde{\mu} \nabla \widetilde{w}^{(\mu)} \Vert_{L^2(B)} + \omega \delta^2 \Vert \widetilde{w}^{(\mu)} \Vert_{L^2(B)}\right).
\end{equation} We need to estimate $\Vert \widetilde{w}^{(\mu)} \Vert_{H^1(B(0,1))}$. Introduce $\widetilde{w}^{(\mu)}_0$ as the solution of
\begin{equation}
\nabla \cdot \left(1+[\sigma_r-1]\mathbf{1}_B + \widetilde{\mu}\right)\nabla \widetilde{w}^{(\mu)}_0 = \nabla \cdot \left( [\sigma_r-1] \mathbf{1}_B \nabla \widetilde{y}\right).
\end{equation} with the condition $\widetilde{w}^{(\mu)}_0(\widetilde{y}) \longrightarrow 0$ as $\widetilde{y} \rightarrow \infty $. Meyers theorem gives:
\begin{equation}\label{estwtilde0}
\Vert \widetilde{w}^{(\mu)}_0 \Vert_{H^1(B(0,1))} \leq C \Vert [\sigma_r -1] \nabla \widetilde{y} \Vert_{L^2(B(0,1))}.
\end{equation}
We can see that $\widetilde{w}^{(\mu)}-\widetilde{w}^{(\mu)}_0$ is a solution of
\begin{equation}
\nabla \cdot \left(1+[\sigma_r-1]\mathbf{1}_B + \widetilde{\mu}\right)\nabla \left(\widetilde{w}^{(\mu)}-\widetilde{w}^{(\mu)}_0\right)  + \omega^2 \delta \left(\widetilde{w}^{(\mu)}-\widetilde{w}^{(\mu)}_0\right) = -\omega^2 \delta \widetilde{w}^{(\mu)}_0.
\end{equation}
We get that
\begin{equation*}
\Vert \widetilde{w}^{(\mu)}-\widetilde{w}^{(\mu)}_0 \Vert_{H^1(B(0,1))} \leq C \omega^2 \delta \Vert \widetilde{w}^{(\mu)}_0 \Vert_{L^2(B(0,1))}.
\end{equation*}
So, using (\ref{estwtilde0}) we get
\begin{equation} \label{estwtilde}
\Vert \widetilde{w}^{(\mu)} \Vert_{H^1(B(0,1))} \leq C \left( 1+ \omega^2 \delta\right).
\end{equation}
Since $\Vert \widetilde{\mu} \nabla\widetilde{w}^{(\mu)}
\Vert_{L^2(B(0,1))} \leq \Vert \widetilde{\mu}\Vert_{L^\infty
(B(0,1))} \Vert \widetilde{w}^{(\mu)}\Vert_{H^1(B(0,1))}$ and
$\Vert \widetilde{\mu}\Vert_{L^\infty(B(0,1))} \leq \Vert \mu
\Vert_\infty$,  using (\ref{estdiffwtild}) and (\ref{estwtilde0})
we get
\begin{equation*}
\Vert \nabla\left(\widetilde{w}^{(\mu)}(\widetilde{y}) -\widetilde{y} -\widetilde{w}(\widetilde{y})\right) \Vert_{L^2(B(0,1))} \leq C \left(\Vert \mu \Vert_\infty + \delta \omega^2(1+\Vert \mu \Vert_\infty + \delta \omega^2)\right),
\end{equation*} which is exactly, as $\Vert \mu \Vert_\infty \rightarrow 0$ and $\delta \rightarrow 0$, for $y\in \Omega_r$
\begin{equation}
\nabla \left( w^{(\mu)}(y)- (y-z_r) \right) = \delta \nabla \widetilde{w}(\frac{y-z_r}{\delta}) + O\left( \delta \Vert \mu \Vert_\infty + (\delta \omega)^2\right).
\end{equation}
\nocite{AmmariIntro}
\nocite{ammari2004reconstruction}
\bibliographystyle{plain}
\bibliography{biblio_final}

\end{document}